\documentclass[11pt,a4paper,final]{article}
\usepackage{amsmath}
\usepackage{amssymb}
\usepackage{amscd}
\usepackage{url}

\usepackage{graphicx}
\usepackage{color}
\usepackage{hyperref}
\hypersetup{bookmarks=true}
\hypersetup{colorlinks=false}

\usepackage[a4paper]{geometry}
\usepackage{amsthm}
\theoremstyle{plain}
\newtheorem{Thm}{Theorem}[section]
\newtheorem{Prop}[Thm]{Proposition}
\newtheorem{Lem}[Thm]{Lemma}
\newtheorem{Cor}[Thm]{Corollary}
\theoremstyle{definition}
\newtheorem{Def}[Thm]{Definition}
\theoremstyle{remark}
\newtheorem{Rmk}[Thm]{Remark}
\newtheorem{Ex}[Thm]{Example}

\newtheorem*{proof*}{proof}

\newcommand{\Sob}[1]{(Sob)$_#1$}
\newcommand{\Dyn}[1]{(Dyn)$_#1$}
\newcommand{\eps}{\varepsilon}

\makeatletter
 
 \@addtoreset{equation}{section}
\makeatother
\usepackage{mathtools}
\mathtoolsset{showonlyrefs=false}

\title{$L^p$-Kato class measures and their relations with Sobolev embedding theorems for Dirichlet spaces}
\author
{%
 Takahiro Mori%
 \thanks
 {%
  Research Institute for Mathematical Sciences,
  Kyoto University, Kyoto, 606-8502, JAPAN.
  \texttt{tmori@kurims.kyoto-u.} \texttt{ac.jp}%
 }
}
\date{}

\begin{document}

\allowdisplaybreaks[4]
\abovedisplayskip=4pt
\belowdisplayskip=4pt

\maketitle

%%%%%%%%%%%%%%%%%%%%%%%%%%%%%%%%%%%%%%%%%%%%%%%%%%%%%%%
%%%Abstract
%%%%%%%%%%%%%%%%%%%%%%%%%%%%%%%%%%%%%%%%%%%%%%%%%%%%%%%
\begin{abstract}
In this paper,
we discuss relationships between
the continuous embeddings of Dirichlet spaces
$(\mathcal{F}, \mathcal{E}_1)$
into Lebesgue spaces and
the integrability of the associated resolvent kernel
$r_\alpha(x, y)$.
For a positive measure $\mu$,
we consider the following two properties;
the first one is that
the Dirichlet space $(\mathcal{F}, \mathcal{E}_1)$ is
continuously embedded into $L^{2p}(E;\mu)$
(which we write as \Sob{p}),
and the second one is that
the family of 1-order resolvent kernels
$\{r_1(x, y)\}_{x\in E}$ is uniformly $p$-th integrable
in $y$ with respect to the measure $\mu$
(which we write as \Dyn{p}).

Under some assumptions,
for a measure $\mu$ satisfying \Dyn{1},
we prove
\Dyn{{p'}} implies \Sob{p}
for $1\leq p \leq p'<\infty$,
and prove
\Sob{{p'}} implies \Dyn{p}
for $1\leq p < p'<\infty$.
To prove these results
we introduce
$L^p$-Kato class,
an $L^p$-version of the set of Kato class measures,
and discuss its properties.
We also give variants of such relations
corresponding to the Gagliardo-Nirenberg type
interpolation inequalities.
As an application,
we discuss the continuity of intersection measures
in time.
\end{abstract}

%%%%%%%%%%%%%%%%%%%%%%%%%%%%%%%%%%%%%%%%%%%%%%%%%%%%%%%
%%%Keywords and Mathematics Subject Classification
%%%%%%%%%%%%%%%%%%%%%%%%%%%%%%%%%%%%%%%%%%%%%%%%%%%%%%%

%31C25    Dirichlet forms
%46E35    Sobolev spaces and other spaces of "smooth'' functions, embedding theorems, trace theorems
%60J35    Transition functions, generators and resolvents
%60J45    Probabilistic potential theory

%%%%%%%%%%%%%%%%%%%%%%%%%%%%%%%%%%%%%%%%%%%%%%%%%%%%%%%
%%%Contents
%%%%%%%%%%%%%%%%%%%%%%%%%%%%%%%%%%%%%%%%%%%%%%%%%%%%%%%
%\setcounter{tocdepth}{2}
%\tableofcontents

%%%%%%%%%%%%%%%%%%%%%%%%%%%%%%%%%%%%%%%%%%%%%%%%%%%%%%%
%%%Main Text
%%%%%%%%%%%%%%%%%%%%%%%%%%%%%%%%%%%%%%%%%%%%%%%%%%%%%%%

%%%%%%%%%%%%%%%%%%%%%%%%%%%%%%%%%%%%%%%%%%%%%%%%%%%%%%%
\section{Introduction}
\label{Sec_intro}
%%%%%%%%%%%%%%%%%%%%%%%%%%%%%%%%%%%%%%%%%%%%%%%%%%%%%%%

In this paper,
we discuss relationships
between
the continuous embeddings of Dirichlet spaces
into Lebesgue spaces
and
the integrability of
the associated resolvent kernel.

The prototype of the relationships
we are focusing on is
the classical Dirichlet integral
$
 \bigl(
  \frac{1}{2}\mathbf{D}
 ,
  H^1(\mathbb{R}^d)
 \bigr)
$
on $\mathbb{R}^d$ and the associated resolvent kernel
$r_\alpha(x, y)$, $x, y\in \mathbb{R}^d$, $\alpha > 0$,
that is,
$H^1(\mathbb{R}^d)$ is the Sobolev space
on $\mathbb{R}^d$,
\begin{equation*}
 \mathbf{D}(u, v)
=
 \sum_{i=1} ^d
 \int_{\mathbb{R}^d}
  \frac{\partial u}{\partial x_i}
  \frac{\partial v}{\partial x_i}
 dx
\quad
 \text{for }
 u, v\in H^1(\mathbb{R}^d)
\end{equation*}

\noindent
and
\begin{equation*}
 r_\alpha(x, y)
=
 \frac{1}{(2\pi)^{d/2}}
 \int_0 ^\infty
  \frac{1}{t^{d/2}}
  \exp
  \biggl\{
  -
  \Bigl(
   \alpha t
  +
   \frac{|x-y|^2}{2t}
  \Bigr)
  \biggr\}
 dt
\quad
 \text{for }
 x, y\in \mathbb{R}^d
,
 \alpha>0
.
\end{equation*}

\noindent
The classical Sobolev embedding theorem
on $\mathbb{R}^d$ is well known:
\begin{equation}
\label{eq_SobRd}
\begin{split}
&
\text{
$
 H^1(\mathbb{R}^d)
$
is continuously embedded into $L^{2p}(\mathbb{R}^d)$
only for
$p\in [1, \infty)$ with
}
\\
&
\text{
$d-p(d-2)\geq 0$
(hence, $2p\leq 2d/(d-2)$ when $d\geq 2$).
}
\end{split}
\end{equation}

\noindent
By an elementary calculation, it holds that
\begin{equation}
\label{eq_resRd}
\text{
$
 \displaystyle
 \sup_{x\in \mathbb{R}^d}
 \int_{\mathbb{R}^d}
  r_1(x, y)^p
 dy
<
 \infty
$
if and only if $d-p(d-2)>0$.
}
\end{equation}
%(see \cite[Theorem A.1]{MR2584458} for example).
Note that $d-p(d-2)$ appears in both conditions
\eqref{eq_SobRd} and \eqref{eq_resRd}.
This means that there is a relation between
the Sobolev embedding theorem and
the integrability of the resolvent.

The main purpose of this paper is to
generalize such relations
from the perspective of the Dirichlet form theory.
Let
$E$ be a locally compact separable metric space,
$m$ be a Radon measure on $E$ with $\mathrm{supp}[m]=E$,
and let
$(\mathcal{E}, \mathcal{F})$ be a Dirichlet form
on $L^2(E; m)$
with the associated resolvent kernel $r_\alpha(x, y)$
with respect to $m$.
We simply write the inner product
$\mathcal{E}_1(u, v) = \mathcal{E}(u, v) + \int_E u v dm$
for $u,v\in \mathcal{F}$.
Suppose $\mu$ is a Radon measure on $E$.
For $p\in [1, \infty)$, we consider two properties;
{
\setlength{\leftmargini}{16mm}
\begin{itemize}

\item[\Sob{p}]
\label{eq_Sobrough}
the Hilbert space
$(\mathcal{F}, \mathcal{E}_1)$ is continuously
embedded into $L^{2p}(E; \mu)$,
that is, there exists a positive constant
$C>0$ such that
$
 \|u\|_{L^{2p}(E;\mu)} ^2
\leq
 C
 \mathcal{E}_1(u, u)
$
for all $u\in \mathcal{F}$,

\item[\Dyn{p}]
\label{eq_Resrough}
it holds that
$
 \displaystyle
 \sup_{x\in E}
 \int_{E}
  r_1(x, y)^p
 \mu(dy)
<
 \infty
$.
\end{itemize}
}

\noindent
The property \Dyn{p} is named after Dynkin,
which is equivalent to $\mathcal{D}^p(X)$
defined later.
The aim of this paper is
to show the following:
under some conditions,
for any measure $\mu$ satisfying \Dyn{1},

{%
\abovedisplayskip=11pt
\belowdisplayskip=11pt

\noindent
\begin{flalign}
\label{eq_Dyn_Sob}
 \hspace{4mm}
 \bullet
 \hspace{1mm}
&
\text{
 if
 \Dyn{{p'}} holds for some $p'\in [1, \infty)$,
 then
 \Sob{p} holds for all $1\leq p \leq p'$,
}
&
\\[3mm]
\label{eq_Sob_Dyn}
 \hspace{4mm}
 \bullet
 \hspace{1mm}
&
\text{
 if
 \Sob{{p'}} holds for some $p'\in [1, \infty)$,
 then
 \Dyn{p} holds for all $1\leq p < p'$.
}
\end{flalign}
}%

\noindent
See Sections \ref{Sec_Sobolev} and \ref{Sec_Sob_pKato}
for precise statements and proofs.
We will prove \eqref{eq_Dyn_Sob} in
Theorem \ref{Thm_Sobolev} (i) and the following
Corollary \ref{Cor_Sobolev}.
We will also prove \eqref{eq_Sob_Dyn} in
Theorem \ref{Thm_Sob_pKato} (i).
Regarding \eqref{eq_Sob_Dyn},
Lemma \ref{Lem_Sob_pKato_2} gives a stronger result
when $\mu$ is the reference measure $m$,
that is,
we obtain the order of decay for the quantity
$
 \sup_{x\in E}
 \int_E
  r_\alpha(x, y)^p
 m(dy)
$
as $\alpha \uparrow \infty$.
By
\eqref{eq_Dyn_Sob} and
\eqref{eq_Sob_Dyn},
\eqref{eq_resRd} follows from
\eqref{eq_SobRd} and
\eqref{eq_SobRd} with $d-p(d-2)>0$
follows from \eqref{eq_resRd}.
In Theorem \ref{Thm_Sobolev} and
Theorem \ref{Thm_Sob_pKato},
we
also give some variants of such relationships.

\vspace{1eM}

When $\mu$ is the reference measure $m$,
the Sobolev inequality has been studied
for various settings;
Euclidean space, Riemannian manifolds,
Lie groups, and so on
(see \cite{MR1872526, MR1386760} for example).
%MR1411441,
%
It is known that the Sobolev inequality is
equivalent to
the ultra-contractivity of
the associated transition semigroup
\cite{MR803094},
the Nash type inequality
\cite{MR898496},
and the capacity isoperimetric inequality
\cite{MR1245225, MR1995492}.

When $p=p'=1$,
our result \eqref{eq_Dyn_Sob} is related to
the theory of the Kato class of measures.
Kato class is introduced to analyse
the Schr\"odinger semigroups
and analyse
integral kernels of semigroups given by
Feynman-Kac functionals
(see \cite{MR644024, MR1185735} for example).
The set of measures satisfying \Dyn{1}
is the so-called the Dynkin class
(1-order version of Green-bounded measures).
The embedding result \eqref{eq_Dyn_Sob}
for $p=p'=1$ is proved
by
Stollmann and Voigt \cite{MR1378151}
via operator theory,
and later
Shiozawa and Takeda \cite{MR2139213}
proved it in terms of Dirichlet forms.

\vspace{1eM}

The organization of the paper is as follows.
In Section \ref{Sec_framework},
we give the framework.
In Section \ref{Sec_pKato},
we introduce an $L^p$-Dynkin class,
which is equivalent to \Dyn{p}, and
an $L^p$-Kato class,
which is an $L^p$-analog of
the classical Kato class.
We will give equivalent conditions of these classes
in terms of heat kernels,
so
we can check that a measure is in the classes
once
an upper bound of  the heat kernel
such as
the (sub-)Gaussian estimate
(see \eqref{eq_subGaussian}) or
the jump type estimate (see \eqref{eq_jump})
holds for a short time.
These estimates are established for many processes.
Regarding the (sub-)Gaussian estimate,
it is obtained for
Brownian motion on a manifold \cite{MR834612},
Brownian motion on a metric measure space
with Riemannian curvature dimension condition
\cite{MR2237207},
Brownian motion on the Sierpi\'nski gasket
\cite{MR966175} and
other diffusions on fractals \cite{MR1668115},
and so on.
Regarding the jump-type estimate,
stable-like processes on $d$ sets \cite{MR2008600}
are studied for example.

In Section \ref{Sec_pdelKato},
we introduce a subclass of the $L^p$-Kato class
(denoted by $\mathcal{K}^{p, \delta}(X)$),
which has additional information
on the order of decay for the quantity
$
 \sup_{x\in E}
 \int_E
  r_\alpha(x, y)^p
 \mu(dy)
$
as $\alpha\uparrow \infty$.
Similarly to
the $L^p$-Dynkin and
the $L^p$-Kato classes,
this subclass
can be characterized via the heat kernel estimate,
and
the above example satisfies the estimate.

In Section \ref{Sec_timechange},
a relation between $L^p$-Dynkin classes
for a process and for its time changed process
is discussed.
This section plays a key role
later to prove \eqref{eq_Sob_Dyn}.
Sections \ref{Sec_Sobolev} and \ref{Sec_Sob_pKato}
are devoted to proving
\eqref{eq_Dyn_Sob} and
\eqref{eq_Sob_Dyn} as mentioned above.
In these sections,
we also prove variants of
\eqref{eq_Dyn_Sob} and \eqref{eq_Sob_Dyn}
corresponding to the $L^p$-Kato class and
the subclass $\mathcal{K}^{p, \delta}(X)$, respectively,
and
the later gives relations with
the Gagliardo-Nirenberg type interpolation inequality.

Section \ref{Sec_ISmeas} is an application
to the intersection of the paths of
independent stochastic processes.
Analysis of the intersection of Brownian paths
was initiated by
Dvoretzky, Erd\H{o}s, Kakutani
\cite{MR0034972, MR0067402}
and
Dvoretzky, Erd\H{o}s, Kakutani and Taylor
\cite{MR0094855}.
They gave
the following dichotomy:
for $p$ independent Brownian motions
$B^{(1)}, \ldots, B^{(p)}$ on $\mathbb{R}^d$,
\begin{equation}
\label{eq_intersection}
\begin{split}
&
\text{
the paths intersect, i.e.,
$
 B^{(1)}(0, \infty)
\cap
 \cdots
\cap
 B^{(p)}(0, \infty)
\not=
 \varnothing
$
almost surely
if
}
\\
&
\text{
$d-p(d-2)>0$,
and
do not intersect almost surely
if
$d-p(d-2)\leq 0$.
}
\end{split}
\end{equation}

\noindent
Note that the same number $d-p(d-2)$ appears
as in \eqref{eq_SobRd} and \eqref{eq_resRd}.

Motivated by problems in statistical physics
such as the configurations of interacting polymers,
a random measure called the intersection local time
has been introduced;
see \cite{MR1229519} for example.
In this paper,
we consider the occupation measure
of the set of intersections
for independent processes $X^{(1)},\ldots, X^{(p)}$
with the same distribution $X$
which is formally written as
\begin{equation}
\notag
%\label{eq_ISmeas}
 \ell^{\mathrm{IS}}_{\boldsymbol{t}}(A)
=
 \int_A
 \biggl[
  \prod_{i=1} ^p
  \int_{0} ^{t_i}
   \delta_x(X^{(i)}(s_i) )
  ds_i
 \biggr]
 m(dx)
\quad
 \text{for }
 A \in \mathcal{B}(E)
\end{equation}
and for
$
 \boldsymbol{t}
=
 (t_1, \ldots, t_p)
\in
 [0, \infty)^p
$,
where
$\delta_x$ is the Dirac measure at $x$,
$m$ is the reference measure of the processes
and
$\mathcal{B}(E)$ is the family of Borel sets in $E$.
We call the measure as
{\it the (mutual) intersection measure}
named after
K\"onig and Mukherjee \cite{MR2999298}.
Here and in the following,
the superscript ``IS'' means ``InterSection''.

In Theorem \ref{Thm_ISconti} we prove the following;
if the reference measure $m$ belongs to
the subclass $\mathcal{K}^{p, \delta}(X)$
introduced in Section \ref{Sec_pdelKato},
then
the measure-valued process
$
 \boldsymbol{t}
\mapsto
 \ell^{\mathrm{IS}}_{\boldsymbol{t}}(dx)
$
has a continuous modification,
and the real-valued process
$
 \boldsymbol{t}
\mapsto
 \langle
  f
 ,
  \ell^{\mathrm{IS}}_{\boldsymbol{t}}
 \rangle
$
has a H\"older continuous modification
for each bounded Borel function $f$.
This is a generalization of
\cite[Section 2.2]{MR2584458},
in which
the results are obtained for
independent Brownian motions.

%\newpage
%\vspace{1eM}

%%%%%%%%%%%%%%%%%%%%%%%%%%%%%%%%%%%%%%%%%%%%%%%%%%%%%%%
\section{$L^p$-Kato class and its variant}
\label{Sec_pKato_pdelKato}
%%%%%%%%%%%%%%%%%%%%%%%%%%%%%%%%%%%%%%%%%%%%%%%%%%%%%%%

In this section,
we give the framework and introduce
the $L^p$-Kato class and its variants.

%%%%%%%%%%%%%%%%%%%%%%%%%%%%%%%%%%%%%%%%%%%%%%%%%%%%%%%
\subsection{Framework}
\label{Sec_framework}
%%%%%%%%%%%%%%%%%%%%%%%%%%%%%%%%%%%%%%%%%%%%%%%%%%%%%%%

Let $E$ be a locally compact, separable metric space
and let $m$ be a Radon measure on $E$
with $\mathrm{supp}[m]=E$.
Let $\partial$ be a point added to $E$ so that
$E_\partial := E \cup \{\partial\}$
is the one-point compactification of $E$.
The point $\partial$ serves
as the cemetery point for $E$.
Suppose $(\mathcal{E}, \mathcal{F})$ is a regular
Dirichlet form on $L^2(E; m)$ and
$X = (\Omega, X_t, \zeta, \mathbb{P}_x)$
is an associated $m$-symmetric Hunt process.
For $\alpha>0$ and $u\in \mathcal{F}$,
we simply write
$
 \mathcal{E}_\alpha(u, u)
=
 \|u\|_{\mathcal{E}_\alpha} ^2
:=
 \mathcal{E}(u, u)
+
 \alpha
 \int_E
  u^2
 dm
$.
In this paper,
%In what follows,
we always take the quasi-continuous version
of the element $u$ of $\mathcal{F}$
(see \cite[Section 2]{MR2778606} for example).

Throughout this paper,
we assume that
the transition kernel $(P_t)_{t>0}$ of $X$
satisfies the absolute continuity condition:
{
\abovedisplayskip=6pt
\belowdisplayskip=6pt
\begin{equation}
\label{eq_AC}
\begin{split}
&
\text{
 $P_t(x, dy)$ is absolutely continuous
 with respect to $m(dy)$ for each
}
\\
&
\text{
 $t>0$ and $x\in E$.
}
\end{split}
\end{equation}
}
%
%Then
%for each $\alpha >0$,
%there exists an $\alpha$-order resolvent kernel
%$r_\alpha(x, y)$ which is defined for all $x, y\in E$
%(see \cite[Lemma 4.2.4]{MR2778606} for example).
%
Note that
the condition \eqref{eq_AC} implies
the measurability of the heat kernel
(see \cite[Theorem 2]{MR960514} for example):

\vspace{-4mm}
{
\abovedisplayskip=6pt
\belowdisplayskip=6pt

\begin{equation}
\label{eq_conti}
\begin{split}
&
 (P_t)_{t>0}
 \text{ admits a heat kernel }
 p_t(x, y)
 \text{ which is jointly measurable on }
\\
&
 (0, \infty)\times E\times E
 \text{ such that }
  p_t(x, y) = p_t(y, x)
  \text{ and }
\\
&
  p_{t+s}(x, y)
 =
  \textstyle
  \int_E
   p_s(x, z)
   p_t(z, y)
  m(dz)
 \text{ for all }
  s, t>0,
  x, y\in E
 .
\end{split}
\end{equation}
}
\vspace{-3mm}

\begin{Rmk}
We may consider a slightly weaker condition
than \eqref{eq_AC}:
{
\abovedisplayskip=6pt
\belowdisplayskip=6pt

\vspace{-1eM}
\begin{equation}
\label{eq_qeAC}
\begin{split}
&
\text{
 There exists a Borel properly exceptional set
 $N$ such that
 $P_t(x, dy)$ is
}
\\
&
\text{
 absolutely
 continuous with respect to $m(dy)$
 for each $t>0$ and $x\in E\setminus N$.
}
\end{split}
\end{equation}
}

\vspace{-3mm}
\noindent
It is known that
the Sobolev inequality implies \eqref{eq_qeAC}
(see \cite[Theorem 4.27]{MR2778606} for example).
Under the condition,
we can obtain the corresponding results of this paper
by replacing
$\sup_{x\in E}$ by
$
 \inf_{\mathrm{Cap}(N)=0}
 \sup_{x\in E\setminus N}
$
as in \cite[(3.3)]{MR1185735},
where $\mathrm{Cap}$ is the ($1$-)capacity of
the Dirichlet form $(\mathcal{E}, \mathcal{F})$.
In this paper,
we do not give detailed calculations
under the assumption \eqref{eq_qeAC}
and
we always impose \eqref{eq_AC} for simplicity.
\end{Rmk}

For each $\alpha >0$,
write the $\alpha$-order resolvent kernel of $X$
by
$
 r_\alpha(x, y)
=
 \int_0 ^\infty
  e^{-\alpha t}
  p_t(x, y)
 dt
$.
We denote by $S_{00}(X)$
the set of positive Borel measures $\mu$
such that
$\mu(E)<\infty$
and
$
 R_1\mu(x)
:=
 \int_E
  r_1(x, y)
 \mu(dy)
$
is uniformly bounded in $x\in E$.
A positive Borel measure $\mu$ on $E$ is
said to be {\it smooth in the strict sense}
if
there exists a sequence $\{E_n\}_{n=1} ^\infty$
of Borel sets increasing to $E$
such that
$
 1_{E_n}\cdot \mu \in S_{00}(X)
$
for each $n$ and
\begin{equation*}
 \mathbb{P}_x
 \Bigl(
  \lim_{n\rightarrow \infty}
  \sigma_{E\setminus E_n}
 \geq
  \zeta
 \Bigr)
=
 1
,
\quad
 \text{for all }x\in E
,
\end{equation*}

\noindent
where
$\sigma_{E\setminus E_n}$ is the first hitting time of
$E\setminus E_n$.
The totality of smooth measures in the strict sense is
denoted by $S_1(X)$.

%\vspace{1eM}
%\newpage

%%%%%%%%%%%%%%%%%%%%%%%%%%%%%%%%%%%%%%%%%%%%%%%%%%%%%%%
\subsection{The class $\mathcal{K}^p(X)$}
\label{Sec_pKato}
%%%%%%%%%%%%%%%%%%%%%%%%%%%%%%%%%%%%%%%%%%%%%%%%%%%%%%%

In this section,
we introduce the $L^p$-version of the Kato class measures.

\begin{Def}
%[$L^p$-Kato class $\mathcal{K}^p(X)$,
%$L^p$-Dynkin class $\mathcal{D}^p(X)$]
\label{Def_pKato}
%\ \par
Let $p\in [1, \infty)$.
For a positive Radon measure $\mu$ on $E$,
$\mu$ is said to be in the
{\it $L^p$-Kato class with respect to $X$}
(in symbols $\mu \in \mathcal{K}^p(X)$)
if
\begin{equation}
\label{eq_pKato}
 \lim_{\alpha\uparrow\infty}
 \sup_{x\in E}
 \int_E
  r_\alpha(x, y)^p
 \mu(dy)
=
 0
\end{equation}

\noindent
and
$\mu$ is said to be in the
{\it $L^p$-Dynkin class with respect to $X$}
(in symbols $\mu \in \mathcal{D}^p(X)$)
if
\begin{equation}
\label{eq_def_pDynkin}
 \sup_{x\in E}
 \int_E
  r_\alpha(x, y)^p
 \mu(dy)
<
 \infty
\quad
 \text{for some }\alpha>0
.
\end{equation}
\end{Def}

\noindent
Clearly $\mathcal{K}^p(X)\subset \mathcal{D}^p(X)$.
The condition \Dyn{p} which we introduced
in Section \ref{Sec_intro}
means nothing else than
the $L^p$-Dynkin class.

%%%%%%%%%%%%%%%%%%%%%%%%%%%%%%%%%%%%%%%%%%%%%%%%%%%%%%%

\begin{Rmk}
\label{Rmk_DefpKato}
\quad
\vspace{-3mm}

\begin{enumerate}
\setlength{\itemsep}{0mm}

\item
$\mathcal{K}^1(X)$ and $\mathcal{D}^1(X)$ are
so-called the set of Kato and Dynkin class measures,
respectively.
The reference measure $m$ always
belongs to $\mathcal{K}^1(X)$.

\item
If $\mu(E)<\infty$, H\"older's inequality gives that
$
 \mu \in \mathcal{K}^{p'}(X)
$
implies
$
 \mu \in \mathcal{K}^{p}(X)
$
for
$
 1\leq p < p'
$.

\item
The $L^p$-Kato and the $L^p$-Dynkin classes are well-defined
up to the quasi-everywhere equivalence of the processes.
Indeed, it holds that
\begin{equation}
\label{eq_sup=qeesssup}
 \sup_{x\in E}
 \int_E
   r_\alpha(x, y)^p
 \mu(dy)
=
 \inf_{N\subset E, \mathrm{Cap}(N)=0}
 \sup_{x\in E\setminus N}
 \int_E
  r_\alpha(x, y)^p
 \mu(dy)
.
\end{equation}

Obviously the left-hand side is equal or grater than
the right-hand side.
To prove the converse inequality,
let $\mu$ be a Radon measure on $E$,
$\alpha > 0$, and $N\subset E$ be a zero-capacity set.
For fixed $x\in E$, $M>0$ and a compact set $K\subset E$,
we write
$
 \nu(dy)
=
 1_K(y) \bigl(r_\alpha(x, y)\wedge M \bigr)^{p-1}
 \mu(dy)
$
and
$
 R_\alpha \nu 
=
 \int_E
  r_\alpha (\cdot, y)
 \nu(dy)
$.
For $z\in E\setminus N$, we have by H\"older's inequality 
\begin{align*}
 R_\alpha \nu(z)
=&
 \int_K
  r_\alpha(z, y)
  \bigl(r_\alpha(x, y)\wedge M \bigr)^{p-1}
 \mu(dy)
\\
\leq&
 \biggl\{
  \sup_{z\in E\setminus N}
  \int_E
   r_\alpha(z, y)^p
  \mu(dy)
 \biggr\}^{\frac{1}{p}}
 \biggl\{
  \int_K
   \bigl(r_\alpha(x, y)\wedge M \bigr)^{p}
  \mu(dy)
 \biggr\}^{\frac{p-1}{p}}
.
\end{align*}

\noindent
Since
$R_\alpha \nu$ is $\alpha$-excessive
and the absolute continuity condition \eqref{eq_AC}
holds, we have

\noindent
\begin{align*}
\notag
 R_\alpha \nu(x)
=&
 \lim_{\varepsilon \downarrow 0}
 e^{-\alpha \varepsilon}
 \mathbb{E}_x [R_\alpha \nu (X_\varepsilon)]
=
 \lim_{\varepsilon \downarrow 0}
 e^{-\alpha \varepsilon}
 \int_{E\setminus N}
  p_\varepsilon(x, z)
  R_\alpha \nu (z)
 m(dz)
\\
\leq&
 \biggl\{
  \sup_{z\in E\setminus N}
  \int_E
   r_\alpha(z, y)^p
  \mu(dy)
 \biggr\}^{\frac{1}{p}}
 \biggl\{
  \int_K
   \bigl(r_\alpha(x, y)\wedge M \bigr)^{p}
  \mu(dy)
 \biggr\}^{\frac{p-1}{p}}
,
\end{align*}

\noindent
which means that
\begin{equation*}
 \int_K
  \bigl(
   r_\alpha(x, y)\wedge M
  \bigr)^{p}
 \mu(dy)
\leq
 \sup_{z\in E\setminus N}
 \int_E
  r_\alpha(z, y)^p
 \mu(dy)
.
\end{equation*}

\noindent
Taking $M\uparrow \infty$ and $K\uparrow E$,
we have the desired inequality
by the dominated convergence theorem.

Hence
one can define the $L^p$-Kato and the $L^p$-Dynkin classes
by the right-hand side of \eqref{eq_sup=qeesssup}
without probabilistic notations and denote them as
$\mathcal{K}^p(\mathcal{E})$ and $\mathcal{D}^p(\mathcal{E})$
in symbols.
Throughout this paper,
we always use the sup-notation for simplicity
because one can find such replacements of $\sup_{x\in E}$ with
$
 \inf_{N\subset E, \mathrm{Cap}(N)=0}
 \sup_{x\in E\setminus N}
$
by similar calculations.

\end{enumerate}
\end{Rmk}

\begin{Ex}
[Brownian motion on $\mathbb{R}^d$]
\label{Ex_BM_pKato}
%\ \par
Suppose $E=\mathbb{R}^d$,
$m$ is the Lebesgue measure on $\mathbb{R}^d$ and
$X$ is a Brownian motion on $\mathbb{R}^d$.
Let $p\in [1, \infty)$ with $d-p(d-2)>0$ and
$\mu$ be a positive Radon measure on $\mathbb{R}^d$.
By the proof of
\cite[Theorem 4.5]{MR644024},
$\mu \in \mathcal{K}^p(X)$
if and only if

\noindent
\begin{align*}
 \lim_{r \downarrow 0}
 \sup_{x\in \mathbb{R}^d}
 \int_{|x-y|<r}
 \frac{\mu(dy)}{|x-y|^{p(d-2)}}
=
 0
,
\quad&
 d\geq 3
,
\\
 \lim_{r \downarrow 0}
 \sup_{x\in \mathbb{R}^d}
 \int_{|x-y|<r}
  \bigl(
   -\log{|x-y|}
  \bigr)^p
 \mu(dy)
=
 0
,
\quad&
 d=2
,
\\
 \sup_{x\in \mathbb{R}^d}
 \int_{|x-y|\leq 1}
 \mu(dy)
<
 \infty
,
\quad&
 d=1
.
\end{align*}
In particular, when $d=1$,
$\mathcal{K}^1(X) = \mathcal{K}^p(X)$
for any $p>1$.

By the above characterization,
we may give a sufficient condition for
$\mathcal{K}^p(X)$.
If
a Borel function $f$ on $\mathbb{R}^d$
satisfies
$
 \sup_{x\in \mathbb{R}^d}
 \int_{|x-y|\leq 1}
  |f(y)|^r
 dy
<
 \infty
$
for some $r > d/(d-p(d-2))$ for $d\geq 2$,
or $r\geq 1$ for $d=1$,
then
the measure
$|f(x)|dx$ is in the class $\mathcal{K}^p(X)$.
This gives an extension of
\cite[Theorem 1.4 (iii)]{MR644024},
in which the result is obtained for $p=1$.
(See also \cite{MR2345907}, in which such results
are obtained under more general heat kernel estimates.)
In particular,
$|x|^{-\beta} dx \in \mathcal{K}^p(X)$
if
$\beta < d-p(d-2)$ for $d\geq 2$,
and $\beta < 1$ for $d=1$.

When $d\geq 3$,
Schechter \cite{MR0447834} introduced
related classes $M_{\alpha, r}$
($\alpha > 0$, $r>1$)
of functions $V$ given by
\begin{equation*}
 \sup_{x\in\mathbb{R}^d}
 \int_{|x-y|\leq 1}
  \frac{|V(y)|^r}{|x-y|^{d-\alpha}}
 dy
<
 \infty
,
\end{equation*}
and \cite{MR644024} studied relations between
$M_{\alpha, r}$
and the classical Kato class $\mathcal{K}^1(X)$.
By H\"older's inequality, we have
$M_{\alpha, r}\subset \mathcal{K}^p(X)$
if $r > \alpha/(d-p(d-2))$.
This is an extension of
\cite[Proposition 4.1, 4.2]{MR644024},
in which the result are obtained for $p=1$.
(Note that
there are typos in \cite{MR644024};
$\beta > 2$ in Proposition 4.1
(resp. $\alpha>2p$ in Proposition 4.2)
should be
$\beta < 2$ (resp. $\alpha < 2p$).)

\end{Ex}

%%%%%%%%%%%%%%%%%%%%%%%%%%%%%%%%%%%%%%%%%%%%%%%%%%%%%%%

\begin{Rmk}
G\"uneysu \cite{10.1093/imrn/rnaa219}
gives another generalization of the Kato class
to show
the H\"older continuity of the Schr\"odinger semigroups.
For $\alpha \in [0, 1]$, a Radon measure $\mu$ on $E$ is said
to be in the $\alpha$-Kato class
in the sense of G\"uneysu
(in symbols $\mu\in \widetilde{\mathcal{K}}^\alpha(X)$) if
\begin{equation*}
 \lim_{t\downarrow 0}
 \sup_{x\in E}
 \int_E
 \biggl(
  \int_0 ^t
   s^{-\frac{\alpha}{2}}
   p_s(x,y) 
  ds
 \biggr)
 \mu(dy)
=
 0
.
\end{equation*}
%区別のためここでは
%$\alpha$-Kato measure in the sense of G\"uneysu
%$\widetilde{\mathcal{K}}^\alpha(X)$
%と書くことにする.

\noindent
Obviously
$\widetilde{\mathcal{K}}^0(X)= \mathcal{K}^1(X)$, 
which is the classical Kato class.
When $X$ is a Brownian motion on $\mathbb{R}^d$,
by the proof of \cite[Theorem 4.5]{MR644024}
$\mu \in \widetilde{\mathcal{K}}^\alpha(X)$
if and only if

\noindent
\begin{align*}
 \lim_{r \downarrow 0}
 \sup_{x\in \mathbb{R}^d}
 \int_{|x-y|<r}
 \frac{\mu(dy)}{|x-y|^{d+\alpha-2}}
=
 0
,
\quad&
 d+\alpha> 2
,
\\
 \lim_{r \downarrow 0}
 \sup_{x\in \mathbb{R}^d}
 \int_{|x-y|<r}
  \bigl(
   -\log{|x-y|}
  \bigr)
 \mu(dy)
=
 0
,
\quad&
 d+\alpha=2
,
\\
 \sup_{x\in \mathbb{R}^d}
 \int_{|x-y|\leq 1}
 \mu(dy)
<
 \infty
,
\quad&
 d+\alpha<2
.
\end{align*}

\noindent
A comparison with Example \ref{Ex_BM_pKato} gives that
the coincidence
$
 \widetilde{\mathcal{K}}^\alpha(X)
=
 \mathcal{K}^p(X)
$
also holds when
$p\geq1$, $d+\alpha-2=p(d-2)>0$,
or when
$p\geq 1$, $d=1$, $0\leq \alpha<1$.
\end{Rmk}

\vspace{1eM}

%%%%%%%%%%%%%%%%%%%%%%%%%%%%%%%%%%%%%%%%%%%%%%%%%%%%%%%

The following proposition
is the $L^p$-version of
\cite[Proposition 3.8]{MR1185735}
in some sense.

\begin{Prop}
\label{Prop_S1}
%\ \par
Let $p\in [1, \infty)$.
It holds that
\begin{equation*}
 \mathcal{D}^p(X)
\subset
 S_1(X)
.
\end{equation*}
\end{Prop}

\begin{proof}
%[Proof of Proposition \ref{Prop_S1}]
%\ \par
Suppose
$\mu \in \mathcal{D}^p(X)$
and
assume first $\mu(E)<\infty$.
As
in Remark \ref{Rmk_DefpKato},
we have $\mu\in \mathcal{D}^1(X)$,
that is,
$
 \sup_{x\in E}
 R_1\mu(x)
<
 \infty
$.
This means $\mu$ is in $S_{00}(X)$
and hence is in $S_1(X)$.

When $\mu\in \mathcal{D}^p(X)$
may not be a finite measure,
take
a sequence $\{E_n\}_{n=1} ^\infty$ of
relatively compact open sets
that go to $E$ as $n\uparrow\infty$.
By the above, we have
$1_{E_n}\cdot \mu \in S_{00}(X)$.
Set
$
 \sigma
=
 \lim_{n\rightarrow \infty}
 \sigma_{E\setminus E_n}
$.
For each $x\in E$,
the quasi-left-continuity of the Hunt process $X$
(see for example, \cite[Appendix A.2]{MR2778606})
implies that
\begin{equation*}
 \lim_{n\rightarrow \infty}
 X_{\sigma_{E\setminus E_n}}
=
 X_{\sigma}
,
\quad
 \text{$\mathbb{P}_x$-a.s. }
 \text{on $\{\sigma <\infty\}$}
\end{equation*}

\noindent
and then
$
 \mathbb{P}_x
 (
  \sigma
 \geq
  \zeta
 )
=
 1
$,
which concludes $\mu \in S_1(X)$.
\end{proof}

\vspace{1eM}

%%%%%%%%%%%%%%%%%%%%%%%%%%%%%%%%%%%%%%%%%%%%%%%%%%%%%%%

The next two propositions characterize
the $L^p$-Dynkin class and
the $L^p$-Kato class
in terms of the heat kernel.

\begin{Prop}
\label{Prop_pDynkin}
%\ \par
Let $p\in [1, \infty)$.
For a Radon measure $\mu$ on $E$,
the following are equivalent:
\vspace{-3mm}

\begin{enumerate}
\setlength{\itemsep}{0mm}
\item[(i)]
 $\mu\in \mathcal{D}^p(X)$,

\item[(ii)]
$
 \displaystyle
 \sup_{x\in E}
 \int_E
  r_\alpha(x, y)^p
 \mu(dy)
<
 \infty
$
for all $\alpha>0$,

\item[(iii)]
$
 \displaystyle
 \sup_{x\in E}
 \int_E
  \biggl(
   \int_0 ^t
    p_s(x, y)
   ds
  \biggr)^p
 \mu(dy)
<
 \infty
$
for some $t>0$,

\item[(iv)]
$
 \displaystyle
 \sup_{x\in E}
 \int_E
  \biggl(
   \int_0 ^t
    p_s(x, y)
   ds
  \biggr)^p
 \mu(dy)
<
 \infty
$
for all $t>0$.
\end{enumerate}
\end{Prop}

%%%%%%%%%%%%%%%%%%%%%%%%%%%%%%%%%%%%%%%%%%%%%%%%%%%%%%%

\begin{proof}
%[Proof of Proposition \ref{Prop_pDynkin}]
%\ \par
Trivially (ii) implies (i) and (iv) implies (iii).

Assume $\mu\in D^p(X)$ and take $\beta > 0$
such that
$
 \sup_{x\in E}
 \int_E
  r_{\beta}(x, y)^p
 \mu(dy)
<
 \infty
$.
For $\alpha>\beta$,
the monotonicity of the resolvent
clearly implies that
$
 \sup_{x\in E}
 \int_E
  r_{\alpha}(x, y)^p
 \mu(dy)
<
 \infty
$.
For $0<\alpha<\beta$,
fix $x\in E$ and set
$
 F(\cdot)
:=
 \int_E
  r_\alpha(x, z)
  r_\beta(z, \cdot)
 m(dz)
$.
Then we have by H\"older's inequality,

\noindent
\begin{align*}
 \int_E
  F(y)^p
 \mu(dy)
=&
 \int_E
  \biggl(
   \int_E
    F(y)^{p-1}
    r_\beta(z, y)
   \mu(dy)
  \biggr)
  r_\alpha(x, z)
 m(dz)
\\
\leq&
 \int_E
  \biggl(
   \int_E
    F(y)^p
   \mu(dy)
  \biggr)^{\frac{p-1}{p}}
  \biggl(
   \int_E
    r_\beta(z, y)^p
   \mu(dy)
  \biggr)^{\frac{1}{p}}
  r_\alpha(x, z)
 m(dz)
\\
\leq&
 \frac{1}{\alpha}
 \biggl(
  \int_E
   F(y)^p
  \mu(dy)
 \biggr)^{\frac{p-1}{p}}
 \biggl(
  \sup_{z\in E}
  \int_E
   r_\beta(z, y)^p
  \mu(dy)
 \biggr)^{\frac{1}{p}}
,
\end{align*}

\noindent
which implies that
\begin{equation*}
 \biggl(
  \int_E
   F(y)^p
  \mu(dy)
 \biggr)^{\frac{1}{p}}
\leq
 \frac{1}{\alpha}
 \biggl(
  \sup_{z\in E}
  \int_E
   r_\beta(z, y)^p
  \mu(dy)
 \biggr)^{\frac{1}{p}}
.
\end{equation*}

\noindent
Hence the resolvent equation
$
 r_\alpha(x, y)
=
 r_\beta(x, y)
+
 (\beta - \alpha)
 \int_E
  r_\alpha(x, z)r_\beta(z, y)
 m(dz)
$
yields that
\begin{align}
\notag
 \biggl(
  \int_E
   r_\alpha(x, y)^p
  \mu(dy)
 \biggr)^{\frac{1}{p}}
\leq&
 \biggl(
  \int_E
   r_\beta(x, y)^p
  \mu(dy)
 \biggr)^{\frac{1}{p}}
+
 (\beta-\alpha)
 \biggl(
  \int_E
   F(y)^p
  \mu(dy)
 \biggr)^{\frac{1}{p}}
\\
\label{eq_pResol}
\leq&
 \frac{\beta}{\alpha}
 \biggl(
  \sup_{z\in E}
  \int_E
   r_\beta(z, y)^p
  \mu(dy)
 \biggr)^{\frac{1}{p}}
<
 \infty
,
\end{align}

\noindent
which concludes (ii).
Moreover, we have
for any $t>0$, $x\in E$ and $\alpha >0$,
\begin{align}
\label{eq_pDynkin_r}
 \int_E
  \biggl(
   \int_0 ^t
    p_s(x, y)
   ds
  \biggr)^p
 \mu(dy)
\leq
 e^{p\alpha t}
 \sup_{x\in E}
 \int_E
  r_\alpha(x, y)^p
 \mu(dy)
,
\end{align}
which concludes (iv).

Next, assume (iii).
Take $t_0>0$ such that
$
 \sup_{x\in E}
 \int_E
  \big(
   \int_0 ^{t_0}
    p_s(x, y)
   ds
  \bigr)^p
 \mu(dx)
<
 \infty
$.
For any $t\leq t_0$ and $a>0$, we have
\begin{equation}
\label{eq_pDynkin}
 \int_E
  \biggl(
   \int_a ^{a+t}
    p_s(x, y)
   ds
  \biggr)^p
 \mu(dy)
\leq
 \sup_{z\in E}
 \int_E
  \biggl(
   \int_0 ^t
    p_s(x, y)
   ds
  \biggr)^p
 \mu(dy)
.
\end{equation}

\noindent
Indeed,
the Chapman-Kolmogorov equation gives
that the left-hand side equals
\begin{align*}
&
 \int_E
  \biggl\{
   \int_0 ^t
   \int_E
    p_a(x, z)
    p_s(z, y)
   m(dz)
   ds
  \biggl\}^p
 \mu(dy)
\\
=&
 \int_E
  \biggl\{
   \int_E
   \biggl(
    \int_0 ^t
     p_s(z, y)
    ds
   \biggl)
   p_a(x, z)
   m(dz)
  \biggl\}^p
 \mu(dy)
.
\end{align*}

\noindent
Applying H\"older's inequality
with the measure $p_a(x, z)m(dz)$,
the above equation is bounded from above by

\noindent
\begin{align*}
&
 \int_E
  \biggl\{
   \int_E
   \biggl(
    \int_0 ^t
     p_s(z, y)
    ds
   \biggl)^p
   p_a(x, z)
   m(dz)
  \biggl\}^{\frac{p}{p}}
  \biggl\{
   \int_E
    1^{\frac{p}{p-1}}
    p_a(x, z)
   m(dz)
  \biggr\}^{\frac{p(p-1)}{p}}
 \mu(dy)
\\
\leq&
 \int_E
  \int_E
   \biggl(
    \int_0 ^t
     p_s(z, y)
    ds
   \biggl)^p
  \mu(dy)
  p_a(x, z)
 m(dz)
\\
\leq&
 \sup_{z\in E}
 \int_E
  \biggl(
   \int_0 ^t
    p_s(z, y)
   ds
  \biggr)^p
 \mu(dy)
,
\end{align*}
where we used $P_a 1 \leq 1$ in the last two lines.
This proves \eqref{eq_pDynkin}.

Now, suppose $t>0$.
By taking large $N$ such that $Nt_0\geq t$,
we have from \eqref{eq_pDynkin}

\noindent
\begin{align*}
 \biggl\{
  \int_E
   \biggl(
    \int_0 ^t
     p_s(x, y)
    ds
   \biggr)^p
  \mu(dy)
 \biggr\}^{\frac{1}{p}}
\leq&
 \sum_{n=0} ^{N-1}
 \biggl\{
  \int_E
   \biggl(
    \int_{nt_0} ^{(n+1)t_0}
     p_s(x, y)
    ds
   \biggr)^p
  \mu(dy)
 \biggr\}^{\frac{1}{p}}
\\
\leq&
 N
 \biggl\{
  \sup_{z\in E}
  \int_E
   \biggl(
    \int_0 ^{t_0}
     p_s(x, y)
    ds
   \biggr)^p
  \mu(dy)
 \biggr\}^{\frac{1}{p}}
,
\end{align*}
which concludes (iv).

Further,
for any $\alpha > 0$, $x\in E$ and $t\leq t_0$,
we have from the triangle inequality,

\noindent
\begin{align}
\notag
 \biggl(
  \int_E
   r_\alpha(x, y)^p
  \mu(dy)
 \biggr)^{\frac{1}{p}}
=&
 \biggl\|
  \int_0 ^\infty
   e^{-\alpha s}
   p_s(x, \cdot)
  ds
 \biggr\|_{L^p(E; \mu)}
\\
\notag
\leq&
 \sum_{n=0} ^\infty
 \biggl\|
  \int_{nt} ^{(n+1)t}
   e^{-\alpha s}
   p_s(x, \cdot)
  ds
 \biggr\|_{L^p(E; \mu)}
\end{align}

\noindent
and by \eqref{eq_pDynkin},
the right-hand side is bounded from above by
\begin{align}
\notag
&
 \sum_{n=0} ^\infty
 e^{-\alpha n t}
 \biggl(
  \sup_{x\in E}
  \int_E
   \biggl(
    \int_0 ^t
     p_s(x, y)
    ds
   \biggr)^p
  \mu(dy)
 \biggr)^{\frac{1}{p}}
\\
\label{eq_pDynkin_p}
=&
 \frac{1}{1-e^{-\alpha t}}
 \biggl(
  \sup_{x\in E}
  \int_E
   \biggl(
    \int_0 ^t
     p_s(x, y)
    ds
   \biggr)^p
  \mu(dy)
 \biggr)^{\frac{1}{p}}
,
\end{align}
which concludes (ii).
\end{proof}

%%%%%%%%%%%%%%%%%%%%%%%%%%%%%%%%%%%%%%%%%%%%%%%%%%%%%%%

\begin{Cor}
\label{Cor_pKato}
%\ \par
Let $p\in [1, \infty)$.
%and suppose that \eqref{eq_conti} holds.
%
Then,
for a Radon measure $\mu$ on $E$,
the following are equivalent:
\vspace{-3mm}

\begin{enumerate}
\setlength{\itemsep}{0mm}

\item[(i)]
$\mu\in \mathcal{K}^p(X)$,

\item[(ii)]
$
 \displaystyle
 \lim_{t\downarrow 0}
 \sup_{x\in E}
 \int_E
  \biggl(
   \int_0 ^t
    p_s(x, y)
   ds
  \biggr)^p
 \mu(dy)
=
 0
$.
\end{enumerate}
\end{Cor}

\begin{proof}
%[Proof of Corollary \ref{Cor_pKato}]
%\ \par
(i) implies (ii)
by letting $t\downarrow 0$
and then $\alpha \uparrow \infty$
in \eqref{eq_pDynkin_r}.
Conversely,
(ii) implies (i)
by letting $\alpha \uparrow \infty$
and then $t\downarrow 0$
in \eqref{eq_pDynkin_p}.
\end{proof}

%%%%%%%%%%%%%%%%%%%%%%%%%%%%%%%%%%%%%%%%%%%%%%%%%%%%%%%

\vspace{1eM}

\begin{Rmk}
\label{Rmk_subproc}
%\ \par
Let $X^{(1)}$ be the 1-subprocess of $X$,
that is,
the $m$-symmetric Markov process
with transition probability
$
 e^{-t}
 p_t(x, y)
 m(dy)
$.
Clearly
$X^{(1)}$ satisfies the absolute continuity condition
\eqref{eq_AC}.
We claim that
$
 \mathcal{K}^{p}(X)
=
 \mathcal{K}^{p}(X^{(1)})
$.
Indeed,
we have the inclusion
$
 \mathcal{K}^{p}(X)
\subset
 \mathcal{K}^{p}(X^{(1)})
$
since
the $\alpha$-order resolvent kernel of $X^{(1)}$ is
$r_{1+\alpha}(x, y)$
and the inequality
$
 r_{1+\alpha}(x, y)
\leq
 r_{\alpha}(x, y)
$
holds.
By applying \eqref{eq_pResol}
with $\beta = \alpha +1$,
we have the converse inclusion
$
 \mathcal{K}^{p}(X^{(1)})
\subset
 \mathcal{K}^{p}(X)
$.
In the same way, we also have
$
 \mathcal{D}^{p}(X)
=
 \mathcal{D}^{p}(X^{(1)})
$.
\end{Rmk}

%\vspace{1eM}
%\newpage

%%%%%%%%%%%%%%%%%%%%%%%%%%%%%%%%%%%%%%%%%%%%%%%%%%%%%%%
\subsection{The class $\mathcal{K}^{p, \delta}(X)$}
\label{Sec_pdelKato}
%%%%%%%%%%%%%%%%%%%%%%%%%%%%%%%%%%%%%%%%%%%%%%%%%%%%%%%

In this section,
we introduce a subclass of $L^p$-Kato class,
which has additional information
on the order of decay of the quantities
$
 \sup_{x\in E}
 \int_E
  r_\alpha(x, y)^p
 \mu(dy)
$
(they have been introduced
in Definition \ref{Def_pKato})
as $\alpha\uparrow \infty$.

\begin{Def}
%\ \par
Let $p\in [1, \infty)$ and $\delta\in (0, 1]$.
For a positive Radon measure $\mu$ on $E$,
$\mu$ is said to be in the
{\it $L^p$-Kato class with order $\delta$}
(in symbols $\mu \in \mathcal{K}^{p, \delta}(X)$)
if
\begin{equation*}
 \sup_{x\in E}
 \biggl(
  \int_E
   r_\alpha(x, y)^p
  \mu(dy)
 \biggr)^{\frac{1}{p}}
=
 O(\alpha^{-\delta})
\quad
 \text{as }
 \alpha \rightarrow \infty
.
\end{equation*}
\end{Def}

\noindent
That is,
there exist constants $C>0$ and $\alpha_0>0$
such that
the left-hand side is bounded from above
by $C \alpha^{-\delta}$
for all $\alpha>\alpha_0$.
Clearly
$
 \mathcal{K}^{p, \delta}(X)
\subset
 \mathcal{K}^p(X)
$.

%%%%%%%%%%%%%%%%%%%%%%%%%%%%%%%%%%%%%%%%%%%%%%%%%%%%%%%

Similarly to
Proposition \ref{Prop_pDynkin} and
Corollary \ref{Cor_pKato},
we can characterize
the set of $L^p$-Kato class measures
with order $\delta$ in terms of the heat kernel.

\begin{Prop}
\label{Prop_pdelKato}
%\ \par
Let $p\in [1, \infty)$ and $\delta\in (0, 1]$.
For a Radon measure $\mu$ on $E$,
the following are equivalent:
\vspace{-3mm}

\begin{enumerate}
\setlength{\itemsep}{0mm}

\item[(i)]
$\mu\in \mathcal{K}^{p, \delta}(X)$,

\item[(ii)]
$
 \displaystyle
 \sup_{x\in E}
 \biggl(
 \int_E
  \biggl(
   \int_0 ^t
    p_s(x, y)
   ds
  \biggr)^p
 \mu(dy)
 \biggl)^\frac{1}{p}
=
 O(t^\delta)
$
as $t\rightarrow 0$.
\end{enumerate}
\end{Prop}

\begin{proof}
%\ \par
By setting $\alpha t = 1$,
(i) implies (ii) from \eqref{eq_pDynkin_r}
and (ii) implies (i) from \eqref{eq_pDynkin_p}.
\end{proof}

%%%%%%%%%%%%%%%%%%%%%%%%%%%%%%%%%%%%%%%%%%%%%%%%%%%%%%%

In the following, we write
\begin{equation}
\label{eq_gamma}
 \gamma(\alpha,\mu,p)
:=
 \sup_{x\in E}
 \biggl(
  \int_E
   r_\alpha(x, y)^p
  \mu(dy)
 \biggr)^{\frac{1}{p}}
\end{equation}

\noindent
and write
\begin{equation}
\label{eq_eta}
 \eta(t,\mu,p)
:=
 \sup_{x\in E}
 \biggl(
  \int_E
   \biggl(
    \int_0 ^t
     p_s(x, y)
    ds
   \biggr)^p
  \mu(dy)
 \biggr)^{\frac{1}{p}}
.
\end{equation}
In the sequel,
we write as $\gamma(\alpha)$ and $\eta(t)$ in short
if there is no danger of confusion.

\begin{Rmk}
\label{Rmk_pdelKato}
\ \par
\vspace{-3mm}

\begin{enumerate}
\setlength{\itemsep}{0mm}
\item
As we see in Corollary \ref{Cor_pKato},
if
$
 \mu
\in
 \mathcal{K}^p(X)
$
then
$
 \eta(t)
<
 \infty
$
for all $t>0$.
Hence,
under \eqref{eq_conti},
$\mu\in \mathcal{K}^{p, \delta}(X)$
if and only if
$
 \sup_{0<t \leq T}
 \bigl\{
  t^{-\delta}
  \eta(t)
 \bigr\}
<
 \infty
$
for some
(also for all) $T>0$.

\item
Let $X^{(1)}$ be the 1-subprocess of $X$.
By the same way
as in Remark \ref{Rmk_subproc},
one can show that
$
 \mathcal{K}^{p, \delta}(X)
=
 \mathcal{K}^{p, \delta}(X^{(1)})
$.

\end{enumerate}
\end{Rmk}

%%%%%%%%%%%%%%%%%%%%%%%%%%%%%%%%%%%%%%%%%%%%%%%%%%%%%%%

\begin{Ex}
[Brownian motion]
\label{Ex_BM_pdelKato}
%\ \par
We continue with Example \ref{Ex_BM_pKato}.
The heat kernel of
the Brownian motion on $\mathbb{R}^d$ is $p_t(x-y)$,
where
\begin{equation*}
 p_t(x)
=
 (2\pi t)^{-\frac{d}{2}}
 e^{-\frac{|x|^2}{2t} }
\quad
 t>0
,
 x\in\mathbb{R}^d
.
\end{equation*}

\noindent
Then,
as in the proof of \cite[Lemma 2.2.4]{MR2584458},
we have for $t>0$,
\begin{equation*}
 \int_{\mathbb{R}^d}
  \biggl(
   \int_0 ^t
    p_s(x)
   ds
  \biggr)^p
 dx
\leq
 (2\pi)^{-\frac{d(p-1)}{2}}
 p^{-\frac{d}{2}}
 \biggl(
  \frac{2p}{2p-d(p-1)}
 \biggr)^p
 t^{\frac{2p-d(p-1)}{2}}
,
\end{equation*}
whenever
$
 2p-d(p-1)
=
  d-p(d-2)
>
 0
$.
Hence,
the Lebesgue measure on $\mathbb{R}^d$ is in
$
 \mathcal{K}^{p, \frac{d-p(d-2)}{2p}}(X)
$
if $d-p(d-2)>0$.
\end{Ex}

\vspace{1eM}

\begin{Ex}
[(Sub-)Gaussian heat kernel estimate]
\label{Ex_subGaussian}
%\ \par
Let $\rho$ denote the metric on $E$, let
$
 \mathrm{diam}(E)
:=
 \sup\{\rho(x, y): x,y\in E\}
$
denote the diameter of $E$,
and denote with
$B(x, r)$ the open ball
with center $x\in E$ and radius $r>0$.
Assume that
$
 \mathrm{diam}(E) = 1
$,
$m(E)<\infty$
and there exist constants
$c_1, c_2>0$ and $d_\mathrm{f}\geq1$ such that
$
 c_1
 r^{d_\mathrm{f}}
\leq
 m(B(x, r))
\leq
 c_2
 r^{d_\mathrm{f}}
$
for all $x\in E$, $r\in (0, 1]$.
We also assume that
$p_t(x, y)$ enjoys
the (sub-)Gaussian heat kernel upper estimate:
there exist constants
$c_3, c_4>0$ and $d_\mathrm{w}\geq 2$ such that

\noindent
\begin{align}
\label{eq_subGaussian}
 p_t(x, y)
\leq
 c_3
 t^{-\frac{d_\mathrm{f}}{d_\mathrm{w}}}
 \exp
 \biggl\{
  -c_4
  \biggl(
   \frac{\rho(x, y)^{\mathrm{d_w}}}{t}
  \biggr)^{\frac{1}{d_\mathrm{w}-1}}
 \biggr\}
\quad
 \text{for all }
 x, y\in E
,
 t\in (0, 1]
.
\end{align}

\noindent
Then, a straightforward calculation gives that
$m\in \mathcal{K}^{p, \delta}(X)$
if
$
 p
 d_\mathrm{w}
 \delta
<
 d_\mathrm{f}
-
 p(d_\mathrm{f} - d_\mathrm{w})
$,
i.e.,
$
 \delta
<
 (d_\mathrm{s} - p(d_\mathrm{s} -2))/{2p}
$
by setting
$
 d_\mathrm{s}
:=
 2d_\mathrm{f}/d_\mathrm{w}
$.
%\begin{equation*}
% \delta
%<
% \frac{d_\mathrm{s} - p(d_\mathrm{s} -2)}{2p}
%\hspace{2mm}
% \biggl(
% =
%  \frac{2p - d_\mathrm{s}(p-1)}{2p}
% \biggr)
%.
%\end{equation*}
%
$d_{\mathrm{f}}$,
$d_{\mathrm{w}}$ and
$d_{\mathrm{s}}$
are the so-called
fractal dimension of $E$
and walk dimension and spectrum dimension
of the process $X$, respectively.
\end{Ex}

\vspace{1eM}

\begin{Ex}
[Jump-type heat kernel estimate]
%\ \par
Under the setting of Example \ref{Ex_subGaussian}
we assume that
$p_t(x, y)$ enjoys
the jump-type heat kernel upper estimate:
there exist constants $c_3>0$ and $d_\mathrm{w}\geq 2$
such that
\begin{align}
\label{eq_jump}
 p_t(x, y)
\leq
 c_3
 \biggl\{
  t^{-\frac{d_\mathrm{f}}{d_\mathrm{w}}}
 \wedge
  \frac{t}{\rho(x, y)^{d_\mathrm{f}+d_\mathrm{w}}}
 \biggr\}
\quad
 \text{for all }
 x, y\in E
,
 t\in (0, 1]
.
\end{align}

\noindent
Then, a straightforward calculation gives that
$m\in \mathcal{K}^{p, \delta}(X)$
if
$
 p
 d_\mathrm{w}
 \delta
<
 d_\mathrm{f}
-
 p(d_\mathrm{f} - d_\mathrm{w})
$,
i.e.,
$
 \delta
<
 (d_\mathrm{s} - p(d_\mathrm{s} -2))/{2p}
$
by setting
$
 d_\mathrm{s}
:=
 2d_\mathrm{f}/d_\mathrm{w}
$.
%\begin{equation*}
% \delta
%<
% \frac{d_\mathrm{s} - p(d_\mathrm{s} -2)}{2p}
%\hspace{2mm}
% \biggl(
% =
%  \frac{2p - d_\mathrm{s}(p-1)}{2p}
% \biggr)
%.
%\end{equation*}
\end{Ex}

\vspace{1eM}

In \cite[Section 1.3]{MR4075013}, these computations
are made in the context of
the intersection measure for such processes.

%\vspace{1eM}

%%%%%%%%%%%%%%%%%%%%%%%%%%%%%%%%%%%%%%%%%%%%%%%%%%%%%%%
\section{$L^p$-Dynkin class with respect to the time changed process}
\label{Sec_timechange}
%%%%%%%%%%%%%%%%%%%%%%%%%%%%%%%%%%%%%%%%%%%%%%%%%%%%%%%

In this section,
we discuss a relation
between $L^p$-Dynkin classes with respect to $X$
and with respect to its time changed process.
The goal of this section is to prove
Proposition \ref{Prop_timechange}.

We first introduce the notation about the time changed
processes of $X$
(for detail, see \cite[Section 6]{MR2778606}
for example).
Let $\{\mathcal{F}_t\}_{t\geq 0}$ is
the minimum completed admissible filtration of $X$,
that is,
$
 \mathcal{F}_t
:=
 \bigcap_{\mu\in \mathcal{P}(E_\partial)}
 \mathcal{F}^\mu_t
$
for $0\leq t< \infty$,
where
$\mathcal{P}(E_\partial)$
is the set of probability measures on $E_\partial$, and
$\mathcal{F}^\mu_t$ is the $\sigma$-algebra
generated by $\{X_s: s\leq t\}$ and
the null sets of the $\mathbb{P}_\mu$-completion of
$\sigma(X_t : 0\leq t < \infty)$.
We also let
$\theta_t$ be the translation operator on $\Omega$,
that is,
$\theta_t$ is a map from $\Omega$ to $\Omega$
such that $X_s\circ\theta_t = X_{s+t}$ for all $s\geq 0$.
A stochastic process $\{A_t\}_{t\geq 0}$
is said to be
a {\it positive continuous additive functional
in the strict sense} (PCAF in abbreviation)
if the following conditions hold:

\vspace{-3mm}

\begin{itemize}
\setlength{\itemsep}{0mm}

\item[(i)]
$A_t(\cdot)$ is $\mathcal{F}_t$-measurable
for all $t \geq 0$,

\item[(ii)]
there exists a set
$
 \Lambda
\in
 \mathcal{F}_\infty
=
 \sigma
 \bigl(
  \bigcup_{t\geq 0}
  \mathcal{F}_t
 \bigr)
$
such that
$\mathbb{P}_x(\Lambda)=1$ for all $x\in E$,
$\theta_t \Lambda \subset \Lambda$ for all $t>0$,
and for each $\omega \in \Lambda$,
$A_{\cdot}(\omega)$ is
a real-valued continuous function satisfying
the following:
$A_0(\omega)=0$,
$A_t(\omega)=A_{\zeta}(\omega)$ for $t\geq \zeta$, and
$A_{t+s}(\omega)=A_t(\omega) + A_s(\theta_t \omega)$
for $t, s\geq0$.
\end{itemize}
\vspace{-3mm}

\noindent
It is known that
there is a one-to-one correspondence between
$S_1(X)$ and the family of PCAF's
(Revuz correspondence)
as follows:
for each $\mu\in S_1(X)$,
there exists a unique PCAF $\{A_t\}_{t\geq 0}$
such that for
any non-negative Borel function $f$ on $E$ and
$\gamma$-excessive function $h$ ($\gamma > 0$),
it holds that
\begin{equation*}
 \int_E
  f(x)
  h(x)
 \mu(dx)
=
 \lim_{t\downarrow 0}
 \frac{1}{t}
 \mathbb{E}_{h\cdot m}
 \biggl[
  \int_0 ^t
   f(X_s)
  dA_s
 \biggr]
,
\end{equation*}
where
$
 \mathbb{E}_{h\cdot m}
 [\hspace{1mm}\cdot\hspace{1mm}]
=
 \int_E
  \mathbb{E}_{x}[\hspace{1mm}\cdot\hspace{1mm}]
  h(x)
 \mu(dx)
$
(see \cite[Theorem 5.1.7]{MR2778606} for example).
We denote by
$A^\mu$ the PCAF corresponding to $\mu\in S_1(X)$.
We
write the fine support of $\mu\in S_1(X)$ as $F$,
that is,
\begin{equation*}
 F
:=
 \{
  x\in E
 :
  \mathbb{P}_x(\tau=0)=1
 \}
,
\quad
 \tau
=
 \inf
 \{
  t>0
 :
  A^\mu_t >0
 \}
.
\end{equation*}

\noindent
%From now on,
%for simplicity we assume that
%the fine support is identical
%to the topological support $\mathrm{supp}[\mu]$,
%and write this as $F$.
%%
%We remark that
%this assumption is not essential,
%and
%the following theorems are valid without it.

For $\mu\in S_1(X)$, denote
$
 \check{X}
=
 (
  \Omega
 ,
  \check{X}_t
 ,
  \check{\zeta}
 ,
  \mathbb{P}_x
 )
$
the time changed process of $X$
with respect to the PCAF $A^\mu$,
that is,

\noindent
\begin{equation*}
 \check{X}_t
=
 X_{\tau_t}
,
\quad
 \tau_t
=
 \inf
 \{
  s>0
 :
  A^\mu_s > t
 \}
,
\quad
 \check{\zeta}
=
 A^\mu_{\zeta}
.
\end{equation*}
Note that
$\check{X}$ is a $\mu$-symmetric Hunt process on $F$.
Write
the $\alpha$-order resolvent of $\check{X}$ by
\begin{equation*}
 \check{R}_\alpha f(x)
=
 \mathbb{E}_x
 \biggl[
  \int_0 ^\infty
   e^{-\alpha t}
   f(\check{X}_t)
  dt
 \biggr]
,
\quad
 f\in \mathcal{B}_b(F)
,
 x\in F
,
\end{equation*}
where $\mathcal{B}_b(F)$ is the set
of bounded Borel functions on $F$.

We note that $\check{X}$ also satisfies
the absolute continuity condition \eqref{eq_AC}:

%%%%%%%%%%%%%%%%%%%%%%%%%%%%%%%%%%%%%%%%%%%%%%%%%%%%%%%

\begin{Lem}
\label{Lem_AC}
Let $\mu$ be a measure in $S_1(X)$.
%whose
%fine support is identical to the topological support.
Then
the time changed process $\check{X}$
satisfies the absolute continuity condition
\eqref{eq_AC}.
\end{Lem}

\begin{proof}
Suppose $\mu(N)=0$.
By
the definition of $S_1(X)$,
we can take
a sequence $\{E_n\}_{n=1} ^\infty$
of Borel sets increasing to $E$
such that
$
 1_{E_n}\cdot \mu \in S_{00}(X)
$
for each $n$ and
\begin{equation}
\label{eq_LemAC1}
 \mathbb{P}_x
 \Bigl(
  \lim_{n\rightarrow \infty}
  \sigma_{E\setminus E_n}
 \geq
  \zeta
 \Bigr)
=
 1
,
\quad
 \text{for all }x\in E
.
\end{equation}

\noindent
For each $n$, the Revuz correspondence
(see \cite[Theorem 5.1.6]{MR2778606} for example)
implies that
\begin{equation*}
 \mathbb{E}_x
 \biggl[
  \int_0 ^\infty
   e^{-\alpha t}
   1_{N\cap E_n}(X_t)
  dA^\mu_t
 \biggr]
=
 R_\alpha
 [1_{N\cap E_n}\cdot \mu]
 (x)
=
 0
\end{equation*}

\noindent
for all $\alpha >0$ and $x\in E$.
By
letting $\alpha \downarrow 0$, we have
\begin{equation*}
 \mathbb{E}_x
 \biggl[
  \int_0 ^\infty
   1_{N\cap E_n}(X_t)
  dA^\mu_t
 \biggr]
=
 0
.
\end{equation*}

\noindent
Since the inclusion
$
 \{
  t< \sigma_{E\setminus E_n}
 \}
\subset
 \bigl\{
  X_t\in E_n \cup \{\partial\}
 \bigr\}
$
holds,
the above equality and \eqref{eq_LemAC1} give that
\begin{align*}
 \mathbb{E}_x
 \biggl[
  \int_0 ^\infty
   1_{N}(X_t)
  dA^\mu_t
 \biggr]
=&
 \lim_{n\rightarrow \infty}
 \mathbb{E}_x
 \biggl[
  \int_0 ^{\sigma_{E\setminus E_n}}
   1_{N}(X_t)
  dA^\mu_t
 \biggr]
\\
\leq&
 \lim_{n\rightarrow \infty}
 \mathbb{E}_x
 \biggl[
  \int_0 ^\infty
   1_{N\cap E_n}(X_t)
  dA^\mu_t
 \biggr]
=
 0
.
\end{align*}

\noindent
Now, we have for every $x\in F$
\begin{align*}
 \check{R}_\alpha 1_N (x)
=&
 \mathbb{E}_x
 \biggl[
  \int_0 ^\infty
   e^{-\alpha t}
   1_N(\check{X}_t)
  dt
 \biggr]
\\
\leq&
 \mathbb{E}_x
 \biggl[
  \int_0 ^\infty
   1_N(\check{X}_t)
  dt
 \biggr]
=
 \mathbb{E}_x
 \biggl[
  \int_0 ^\infty
   1_N(X_t)
  dA^\mu_t
 \biggr]
=
 0
,
\end{align*}
which concludes the absolute continuity condition
for $\check{X}$
(see \cite[Theorem 4.2.4]{MR2778606} for example).
\end{proof}

\vspace{1eM}

%%%%%%%%%%%%%%%%%%%%%%%%%%%%%%%%%%%%%%%%%%%%%%%%%%%%%%%

The next proposition will play a key role later
to prove Theorem \ref{Thm_Sob_pKato} (i),
one of our main results.
Roughly it means that,
if $\mu$ is $L^p$-Dynkin with respect to
the time changed process,
then
$\mu$ is $L^p$-Dynkin with respect to
the original process.

\begin{Prop}
\label{Prop_timechange}
%\ \par
Let $p\in (1, \infty)$ and
let $\mu$ be a measure in $S_1(X)$.
%whose
%fine support is identical to the topological support.
Then,
$\mu \in \mathcal{D}^1(X)$ and
$\mu\in \mathcal{D}^p(\check{X})$
imply
$\mu \in \mathcal{D}^p(X)$.
\end{Prop}

To prove this,
we first consider the transient version
of Proposition \ref{Prop_timechange}.

\begin{Def}
%\ \par
When $(\mathcal{E}, \mathcal{F})$
is transient,
$\mu$ is said to be {\it Green-bounded}
(in symbols $\mu\in \mathcal{D}_0(X)$) if
\begin{equation}
\label{eq_pGreenbdd}
 \sup_{x\in E}
 \int_E
  r_0(x, y)
 \mu(dy)
<
 \infty
,
\end{equation}
where
$
 r_0(x, y)
=
 \lim_{\alpha \downarrow 0}
 r_\alpha(x, y)
$.
\end{Def}

\begin{Lem}
\label{Lem_timechange}
%\ \par
Suppose $(\mathcal{E}, \mathcal{F})$ is transient.
Let
$p\in (1, \infty)$ and $\mu$ be a smooth measure
in $S_1(X)$.
%whose fine support is identical
%to the topological support.
Then,
$\mu \in \mathcal{D}_0(X)$
and
$\mu \in \mathcal{D}^p(\check{X})$
imply
$\mu \in \mathcal{D}^p(X)$.
\end{Lem}

%%%%%%%%%%%%%%%%%%%%%%%%%%%%%%%%%%%%%%%%%%%%%%%%%%%%%%%

\begin{proof}
[Proof of Lemma \ref{Lem_timechange}]
%\ \par
Denote $\check{r}_\alpha(x, y)$ the $\alpha$-order
resolvent kernel of $\check{X}$.
First,
we claim that for all $x\in F$,
\begin{equation}
\label{eq_check}
 \check{r}_0(x, y)
=
 r_0(x, y)
\quad
 \text{for $\mu$-a.e. $y\in F$}.
\end{equation}

\noindent
Take a sequence $\{E_n\}_{n=1} ^\infty$
of Borel sets increasing to $E$
as in the definition of $\mu\in S_1(X)$.
Let $f$ be a non-negative Borel function on $F$.
By a similar argument as in the proof of
Lemma \ref{Lem_AC}, we have
\begin{align*}
 \mathbb{E}_x
 \biggl[
  \int_0 ^\infty
   1_{\bigcup_{n=1} ^\infty E_n}(X_t)
   f(X_t)
  dA^\mu_t
 \biggr]
=
 R_0[ f\cdot \mu] (x)
\end{align*}
for all $x\in E$.
The left-hand side of the above equation
is equal to
$
 \mathbb{E}_x
 \bigl[
  \int_0 ^\infty
   f(X_t)
  dA^\mu_t
 \bigr]
$.
Indeed,
by the inclusion
$
 \{t< \sigma_{E\setminus E_n}\}
\subset
 \bigl\{
  X_t
 \in
  E_n \cup \{\partial\}
 \bigr\}
$
we have

\noindent
\begin{align*}
 \mathbb{E}_x
 \biggl[
  \int_0 ^\infty
   f(X_t)
  dA^\mu_t
 \biggr]
=&
 \lim_{n\rightarrow \infty}
 \mathbb{E}_x
 \biggl[
  \int_0 ^{\sigma_{E\setminus E_n}}
   f(X_t)
  dA^\mu_t
 \biggr]
\\
\leq&
 \mathbb{E}_x
 \biggl[
  \int_0 ^\infty
   1_{\bigcup_{n=1} ^\infty E_n}(X_t)
   f(X_t)
  dA^\mu_t
 \biggr]
.
\end{align*}
The converse inequality is trivial.
Hence,
\eqref{eq_check} follows from
the definition of $0$-order resolvent
\begin{equation*}
 R_0[f\cdot \mu](x)
=
 \int_F
  r_0(x, y)
  f(y)
 \mu(dy)
\end{equation*}
for $x\in E$
and from the equality
\begin{equation*}
 \mathbb{E}_x
 \biggl[
  \int_0 ^\infty
   f(X_t)
  dA^\mu_t
 \biggr]
=
 \mathbb{E}_x
 \biggl[
  \int_0 ^\infty
   f(X^\mu_t)
  dt
 \biggr]
=
 \int_F
  \check{r}_0(x, y)
 \mu(dy)
\end{equation*}
for $x\in F$
which is obtained from the change of variables.

Next, we recall the resolvent equation
\begin{align*}
 \check{r}_0(x, y)
=
 \check{r}_\alpha(x, y)
+
 \alpha
 \int_F
  \check{r}_0(x, z)
  \check{r}_\alpha(z, y)
 \mu(dz)
,
\quad
 \alpha >0
,
 x, y\in F
.
\end{align*}

\noindent
A similar calculation as \eqref{eq_pResol}
gives that, for $x\in F$

\noindent
\begin{align*}
\begin{split}
&
 \biggl(
  \int_F
   \check{r}_0(x, y)^p
  \mu(dy)
 \biggr)^{\frac{1}{p}}
\\
\leq&
 \biggl(
  1
 +
  \alpha
  \sup_{x\in F}
  \int_F
   \check{r}_0(x, y)
  \mu(dy)
 \biggr)
 \biggl(
  \sup_{x\in F}
  \int_{F}
   \check{r}_\alpha(x, y)^p
  \mu(dy)
 \biggr)^{\frac{1}{p}}
.
\end{split}
\end{align*}

\noindent
By combining this with \eqref{eq_check},
we have
\begin{align}
\label{eq_check_bdd}
\begin{split}
&
 \biggl(
  \sup_{x\in F}
  \int_F
   r_0(x, y)^p
  \mu(dy)
 \biggr)^{\frac{1}{p}}
\\
\leq&
 \biggl(
  1
 +
  \alpha
  \sup_{x\in E}
  \int_E
   r_0(x, y)
  \mu(dy)
 \biggr)
 \biggl(
  \sup_{x\in F}
  \int_{F}
   \check{r}_\alpha(x, y)^p
  \mu(dy)
 \biggr)^{\frac{1}{p}}
.
\end{split}
\end{align}

\noindent
The right-hand side of \eqref{eq_check_bdd} is
finite because of the assumptions
$\mu\in \mathcal{D}_0(X)$
and
$\mu\in \mathcal{D}^p(\check{X})$.

Now we will show that
\begin{equation}
\label{eq_EbddF}
 \sup_{x\in E}
 \int_F
  r_0(x, y)^p
 \mu(dy)
=
 \sup_{z\in F}
 \int_{F}
  r_0(z, y)^p
 \mu(dy)
.
\end{equation}

\noindent
Let $x\in F$.
Define the $0$-order hitting distribution
$H^0_F(x, dz)$ by
\begin{equation*}
 H^0_F(x, A)
=
 \mathbb{E}_x
 [
  1_A(X_{\sigma_F})
 ;
  \sigma_F <\infty
 ]
,
\quad
 \text{for }x\in E,
 \hspace{1ex}
 A\in\mathcal{B}(E)
.
\end{equation*}

\noindent
Then we can see that
\begin{equation*}
 r_0(x, y)
=
 \int_{F}
  r_0(z, y)
 H^0_F(x, dz)
\quad
 \text{q.e. } y\in E
.
\end{equation*}
Indeed,
for a non-negative Borel function $f$ on $E$ we have
from the strong Markov property,

\noindent
\begin{align*}
&
 \int_E
  r_0(x, y)
  f(y)
 m(dy)
\\
=&
 \mathbb{E}_x
 \biggl[
  \int_0 ^\infty
   f(X_t)
  dt
 \biggr]
=
 \mathbb{E}_x
 \biggl[
  \int_{\sigma_F} ^\infty
   f(X_t)
  dt
 \biggr]
\\
=&
 \mathbb{E}_x
 \bigl[
  R_0 f (X_{\sigma_F})
 ;
  \sigma_F < \infty
 \bigr]
=
 \int_E
  \biggl(
   \int_{F}
    r_0(z, y)
   H^0_F(x, dz)
  \biggr)
  f(y)
 m(dy)
\end{align*}

\noindent
and hence the equality holds for $m$-a.e. $x\in E$.
The desired equality for $\text{q.e. }x\in E$
follows from the fact that
the functions on both hand sides are
$0$-excessive in $y\in E$.
By
applying H\"older's inequality to the measure
$H^0_F(x, dz)$, we have

\noindent
\begin{align*}
 \int_E
  r_0(x, y)^p
 \mu(dy)
=&
 \int_F
  r_0(x, y)^{p-1}
  \biggl(
   \int_{F}
    r_0(z, y)
   H^0_F(x, dz)
  \biggr)
 \mu(dy)
\\
=&
 \int_{F}
 \int_F
  r_0(x, y)^{p-1}
  r_0(z, y)
 \mu(dy)
 H^0_F(x, dz)
\\
\leq&
 \int_{F}
  \biggl(
   \int_F
    r_0(x, y)^p
   \mu(dy)
  \biggr)^{\frac{p-1}{p}}
  \biggl(
   \int_F
    r_0(z, y)^p
   \mu(dy)
  \biggr)^{\frac{1}{p}}
 H^0_F(x, dz)
\end{align*}

\noindent
and the right-hand side is bounded from above by
\begin{align*}
 \biggl(
  \int_F
   r_0(x, y)^p
  \mu(dy)
 \biggr)^{\frac{p-1}{p}}
 \biggl(
  \sup_{z\in F}
  \int_F
   r_0(z, y)^p
  \mu(dy)
 \biggr)^{\frac{1}{p}}
\end{align*}
because of
$
 H^0_F(x, F)
\leq
 1
$.
Hence
we obtain \eqref{eq_EbddF}.

By combining
\eqref{eq_check_bdd} with \eqref{eq_EbddF}, we have
\begin{align*}
&
 \biggl(
  \sup_{x\in E}
  \int_F
   r_\alpha(x, y)^p
  \mu(dy)
 \biggr)^{\frac{1}{p}}
\\
\leq&
 \biggl(
  \sup_{x\in E}
  \int_F
   r_0(x, y)^p
  \mu(dy)
 \biggr)^{\frac{1}{p}}
\\
\leq&
 \biggl(
  1
 +
  \alpha
  \sup_{x\in F}
  \int_F
   r_0(x, y)
  \mu(dy)
 \biggr)
 \biggl(
  \sup_{x\in F}
  \int_{F}
   \check{r}_\alpha(x, y)^p
  \mu(dy)
 \biggr)^{\frac{1}{p}}
,
\end{align*}
which completes the proof.
\end{proof}

\vspace{1eM}

%%%%%%%%%%%%%%%%%%%%%%%%%%%%%%%%%%%%%%%%%%%%%%%%%%%%%%%

We now prove
Proposition \ref{Prop_timechange}.

\begin{proof}
[Proof of Proposition \ref{Prop_timechange}]
%\ \par
Let
$Y=X^{(1)}$ be the 1-subprocess of $X$
defined in Remark \ref{Rmk_subproc}.
We can find that
$\mu$ is in $S_1(Y)$.
%and
%its fine support with respect to $Y$ is identical
%to the topological support.
%
We can also find that
the assumptions
$\mu \in \mathcal{D}(X)$ and
$\mu\in \mathcal{D}^p(\check{X})$
imply
$\mu \in \mathcal{D}_0(Y)$ and
$\mu\in \mathcal{D}^p(\check{Y})$.
Since
$Y$ is transient,
Lemma \ref{Lem_timechange} gives that
$\mu\in \mathcal{D}^p(Y)$.
The conclusion
follows from the equality
$
 \mathcal{D}^p(X)
=
 \mathcal{D}^p(Y)
$,
which has already been noted in Remark
\ref{Rmk_subproc}.
\end{proof}

%\newpage
%\vspace{1eM}

%%%%%%%%%%%%%%%%%%%%%%%%%%%%%%%%%%%%%%%%%%%%%%%%%%%%%%%
\section{Main results}
\label{Sec_result}
%%%%%%%%%%%%%%%%%%%%%%%%%%%%%%%%%%%%%%%%%%%%%%%%%%%%%%%

%This is the main part of this paper.
In this section, we give relations between
the $L^p$-Dynkin classes and the Sobolev embeddings
as we introduced in
\eqref{eq_Dyn_Sob} and \eqref{eq_Sob_Dyn}.
We also give variants of such relations
corresponding to the $L^p$-Kato class and
the $L^p$-Kato class with order $\delta$,
respectively.

%%%%%%%%%%%%%%%%%%%%%%%%%%%%%%%%%%%%%%%%%%%%%%%%%%%%%%%
\subsection{$L^p$-Kato implies the Sobolev embedding}
\label{Sec_Sobolev}
%%%%%%%%%%%%%%%%%%%%%%%%%%%%%%%%%%%%%%%%%%%%%%%%%%%%%%%

In this section,
we first give
an $L^p$-version of
the Stollmann-Voigt inequality
and
some variants in Theorem \ref{Thm_Sobolev}.
As a consequence,
Corollary \ref{Cor_Sobolev} proves
the assertion \eqref{eq_Dyn_Sob}
introduced in Section \ref{Sec_intro},
that is, for
a measure $\mu\in \mathcal{D}^1(X)$ and
$1\leq p\leq p'$,
$\mu\in \mathcal{D}^{p'}(X)$ implies that
$(\mathcal{F}, \mathcal{E}_1)$ is continuously
embedded into $L^{2p}(E;\mu)$.
We also give a Rellich-Kondrachov type
compact embedding theorem
(Corollary \ref{Cor_cptemb}).

%%%%%%%%%%%%%%%%%%%%%%%%%%%%%%%%%%%%%%%%%%%%%%%%%%%%%%%

\begin{Thm}
\label{Thm_Sobolev}
%\ \par
Let $p\in [1, \infty)$ and $\mu \in \mathcal{D}^p(X)$.

\vspace{-3mm}

\begin{itemize}
\setlength{\itemsep}{0mm}

\item[(i)]
It holds that
\begin{equation}
\label{eq_pSV}
 \|
  u
 \|_{L^{2p}(E;\mu)} ^2
\leq
 \biggl(
  \sup_{x\in E}
  \int_E
   r_\alpha(x, y)^p
  \mu(dy)
 \biggr)^{\frac{1}{p}}
 \mathcal{E}_\alpha(u, u)
\end{equation}
for any $u\in \mathcal{F}$ and $\alpha > 0$.
In particular,
the Hilbert space
$(\mathcal{F}, \mathcal{E}_1)$ is
continuously embedded into $L^{2p}(E;\mu)$.

\item[(ii)]
If $\mu \in \mathcal{K}^p(X)$,
then it holds that
\begin{equation}
\label{eq_Kato_Sob}
 \|
  u
 \|_{L^{2p}(E;\mu)} ^2
\leq
 \eps
 \mathcal{E}_1(u, u)
+
 K(\eps)
 \|
  u
 \|_{L^{2}(E;m)} ^2
\end{equation}
for any $u\in \mathcal{F}$ and $\eps > 0$,
where
$K$ is a positive function on $(0, \infty)$
satisfying
$\eps^{-1} K(\eps) \uparrow \infty$
as $\eps \downarrow 0$.

\item[(iii)]
If $\mu \in \mathcal{K}^{p, \theta}(X)$
for some $\theta\in (0, 1]$,
then \eqref{eq_Kato_Sob} holds for
a function
$K(\eps) = A\eps^{-\frac{1-\theta}{\theta}}$,
where $A$ is a positive constant.

In particular, it holds that
\begin{equation}
\label{eq_thetaKato_Sob}
 \|
  u
 \|_{L^{2p}(E;\mu)}
\leq
 B
 \sqrt{\mathcal{E}_1(u, u)}^{(1-\theta)}
 \|
  u
 \|_{L^{2}(E;m)} ^{\theta}
\end{equation}
for any $u\in \mathcal{F}$,
where
$B$ is an another positive constant.
\end{itemize}

\end{Thm}

Note that
\eqref{eq_Kato_Sob} is similar to the notion of
compactly boundedness in \cite{MR369884},
and that
\eqref{eq_thetaKato_Sob} is
so-called the interpolation type inequality.

%%%%%%%%%%%%%%%%%%%%%%%%%%%%%%%%%%%%%%%%%%%%%%%%%%%%%%%

\begin{proof}
[Proof of Theorem \ref{Thm_Sobolev}]
%\ \par
We first assume $\mu\in \mathcal{D}^p(X)$
and prove \eqref{eq_pSV}.
The case $p=1$ is exactly
the Stollmann-Voigt inequality
(see \cite{MR2139213, MR1378151}),
%\cite[Exercise 6.4.4]{MR2778606},
so we assume $p>1$.
By
the regularity of $(\mathcal{E}, \mathcal{F})$,
it suffices to prove \eqref{eq_pSV} for
$u\in \mathcal{F}\cap C_0(E)$,
where
$C_0(E)$ is the set of continuous functions on $E$
with compact support.
Fix
$u\in \mathcal{F}\cap C_0(E)$
and $\alpha >0$.
Define
a finite measure $\nu$ on $E$ by
$
 \nu(dy)
=
 u^{2p-2}(y)\mu(dy)
$.
By H\"older's inequality,
we have for $x\in E$,
\begin{align*}
 R_\alpha\nu(x)
=&
 \int_E
  r_\alpha(x, y)
  u^{2p-2}(y)
 \mu(dy)
\leq
 \biggl(
  \sup_{x\in E}
  \int_E
   r_\alpha(x, y)^p
  \mu(dy)
 \biggr)^{\frac{1}{p}}
 \biggl(
  \int_E
   u^{2p}
  d\mu
 \biggr)^{\frac{p-1}{p}}
.
\end{align*}

\noindent
By applying
the inequality \eqref{eq_pSV} with $p=1$
and $\nu \in S_1(X)$,
we have
\begin{align*}
 \int_E
  u^{2p}
 d\mu
=&
 \int_E
  u^2
 d\nu
\\
\leq&
 \|
  R_\alpha \nu
 \|_\infty
 \mathcal{E}_\alpha(u, u)
\\
\leq&
 \biggl(
  \sup_{x\in E}
  \int_E
   r_\alpha(x, y)^p
  \mu(dy)
 \biggr)^{\frac{1}{p}}
 \biggl(
  \int_E
   u^{2p}
  d\mu
 \biggr)^{\frac{p-1}{p}}
 \mathcal{E}_\alpha(u, u)
,
\end{align*}
which concludes \eqref{eq_pSV}.

\vspace{1eM}

We next assume $\mu\in \mathcal{K}^p(X)$
and prove \eqref{eq_Kato_Sob}.
Recall
the notation $\gamma(\alpha)$ introduced in
\eqref{eq_gamma}.
We may assume $\gamma(\alpha)>0$ for all $\alpha>0$.
Indeed,
if $\gamma(\alpha_0)=0$ for some $\alpha_0>0$,
then
$R_\alpha \mu(x) = 0$ for all $\alpha \geq \alpha_0$
and $x\in E$, and hence $\mu=0$.
In this case
\eqref{eq_Kato_Sob} clearly holds for
$K(\eps)\equiv 1$.

Note that the function $\gamma(\cdot)$
is continuous.
Indeed,
for $0<\alpha_0<\alpha < \beta$,
the monotonicity of $\gamma(\cdot)$
and \eqref{eq_pResol} imply that
\begin{align*}
 0
\leq
 \gamma(\alpha) - \gamma(\beta)
\leq
 \frac{\beta}{\alpha}
 \gamma(\beta)
-
 \gamma(\beta)
=
 \frac{\gamma(\beta)}{\alpha}(\beta - \alpha)
\leq
 \frac{\gamma(\alpha_0)}{\alpha_0}(\beta - \alpha)
.
\end{align*}

Define the right continuous inverse $\gamma^{-1}$
of $\gamma$ by
\begin{equation*}
 \gamma^{-1}(\eps)
:=
 \inf
 \{
  \alpha>0
 :
  \gamma(\alpha) \leq \eps
 \},
\quad
 \eps > 0
.
\end{equation*}
We have
$\gamma^{-1}(\eps)\uparrow \infty$
as $\eps \downarrow 0$
because of the assumption $\mu\in \mathcal{K}^p(X)$.
We also have
$\gamma(\gamma^{-1}(\eps)) \leq \eps$
because of the right continuity of $\gamma^{-1}$.
By
substituting $\alpha = \gamma^{-1}(\eps)$ in
\eqref{eq_pSV},
we obtain \eqref{eq_Kato_Sob} with
$K(\eps)= \eps \gamma^{-1}(\eps)$.

\vspace{1eM}

Finally,
we assume $\mu\in \mathcal{K}^{p, \theta}(X)$
for $\theta\in (0, 1]$
and prove \eqref{eq_thetaKato_Sob}.
By
the definition of $\mathcal{K}^{p, \theta}(X)$,
there exists a constant $C>0$ such that
$
 \gamma(\alpha+1)
\leq
 C(\alpha+1)^{-\theta}
\leq
 C\alpha^{-\theta}
$
for all $\alpha>0$.
Then,
for $u\in \mathcal{F}$
we have from \eqref{eq_pSV}
\begin{align*}
 \|
  u
 \|_{L^{2p}(E;\mu)} ^2
\leq&
 \gamma(\alpha + 1)
 \mathcal{E}_{\alpha+1}(u, u)
\leq
 C\alpha^{-\theta}
 \mathcal{E}_{1}(u, u)
+
 C\alpha^{1-\theta}
 \|
  u
 \|_{L^{2}(E;m)} ^2
.
\end{align*}
By substituting
$
 \alpha
=
 C^{\frac{1}{\theta}}
 \eps^{-\frac{1}{\theta}}
$
in the right-hand side of the above inequality,
we obtain \eqref{eq_Kato_Sob} with
$
 K(\eps)
=
 C^{\frac{1}{\theta}}
 \eps^{-\frac{1-\theta}{\theta}}
$.
The interpolation inequality
\eqref{eq_thetaKato_Sob} follows
by taking infimum over $\eps>0$.
\end{proof}

%%%%%%%%%%%%%%%%%%%%%%%%%%%%%%%%%%%%%%%%%%%%%%%%%%%%%%%

By combining Theorem \ref{Thm_Sobolev} (i)
with H\"older's inequality, we have the following
corollary.

\begin{Cor}
\label{Cor_Sobolev}
%\ \par
Let $p'\in [1, \infty)$.
Then,
for any measure
$
 \mu
\in
 \mathcal{D}^1(X)\cap\mathcal{D}^{p'}(X)
$
the Hilbert space
$(\mathcal{F}, \mathcal{E}_1)$ is continuously embedded into
$L^{2p}(E;\mu)$
for all $1\leq p\leq p'$.
\end{Cor}

\vspace{1eM}
%\newpage

%%%%%%%%%%%%%%%%%%%%%%%%%%%%%%%%%%%%%%%%%%%%%%%%%%%%%%%

At the end of this section,
we give
a Rellich-Kondrachov type compact embedding theorem.
The following corollary is a generalization of
\cite[Corollary 4.5]{MR4009422},
where the statement is proved for $p=1$ and $\mu=m$.

\begin{Cor}
\label{Cor_cptemb}
%\ \par
%In addition to
%the assumption in Proposition \ref{Prop_Sobolev},
Assume $X$ satisfies
\vspace{-3mm}

\begin{itemize}
\setlength{\itemsep}{0mm}

\item
(resolvent strong Feller property)
$
 R_1(\mathcal{B}_b(E))
\subset
 C_b(E)
$,
where $C_b(E)$ is
the set of bounded continuous functions on $E$,
and

\item
(tightness)
for any $\varepsilon >0$,
there exists a compact set $K\subset E$
such that
$
 \sup_{x\in E}
 R_1 1_{K^c}(x)
<
 \varepsilon
$.
\end{itemize}
\vspace{-3mm}

\noindent
Let $p\in [1, \infty)$
and suppose $\mu \in \mathcal{K}^p(X)$.
Then
the Hilbert space $(\mathcal{F}, \mathcal{E}_1)$
is compactly embedded into
$L^2(E;m)$ and $L^{2p}(E;\mu)$.
\end{Cor}

\begin{proof}
%\ \par
Suppose
$
 \{u_n\}_{n=1} ^\infty
\subset
 \mathcal{F}
$
is bounded in $(\mathcal{F}, \mathcal{E}_1)$.
By
\cite[Corollary 4.5]{MR4009422},
$
 (\mathcal{F}, \mathcal{E}_1)
$
is compactly embedded into $L^2(E; m)$.
(We remark that
the irreducibility assumption is not needed
to prove \cite[Corollary 4.5]{MR4009422}.)
Hence
we can take $u\in L^2(E; m)$
and a subsequence $\{u_{n(k)}\}_{k=1} ^\infty$
such that
$u_{n(k)}$ converges to $u$ in $L^2(E; m)$
as $k\rightarrow \infty$.

Recall the notation $\gamma(\alpha)$
introduced in \eqref{eq_gamma}.
By \eqref{eq_pSV}, we have
\begin{align*}
&
 \|
  u_{n(k)}
 -
  u_{n(l)}
 \|_{L^{2p}(E; \mu)}^2
\\
\leq&
 \gamma(\alpha)
 \|
  u_{n(k)}
 -
  u_{n(l)}
 \|_{\mathcal{E}_\alpha}^2
\\
=&
 \gamma(\alpha)
 \mathcal{E}
 (
  u_{n(k)} - u_{n(l)}
 ,
  u_{n(k)} - u_{n(l)}
 )
+
 \alpha
 \gamma(\alpha)
 \|
  u_{n(k)}
 -
  u_{n(l)}
 \|_{L^2(E;m)}^2
\\
\leq&
 4\gamma(\alpha)
 \sup_n
 \|
  u_n
 \|_{\mathcal{E}_1}^2
+
 \alpha
 \gamma(\alpha)
 \|
  u_{n(k)}
 -
  u_{n(l)}
 \|_{L^2(E;m)}^2
.
\end{align*}

\noindent
By letting $k, l\rightarrow \infty$,
and then $\alpha \rightarrow \infty$,
we find that $\{u_{n(k)}\}$ is
Cauchy in $L^{2p}(E;\mu)$,
and this completes the proof.
\end{proof}

\vspace{1eM}

%%%%%%%%%%%%%%%%%%%%%%%%%%%%%%%%%%%%%%%%%%%%%%%%%%%%%%%

\begin{Rmk}
%\ \par
By Corollary \ref{Cor_pKato},
assumption (A5) of \cite{MR4075013}
is equivalent to
the reference measure $m$ belonging to
the $L^p$-Kato class.
We can also see that
assumption (A4) of \cite{MR4075013} is used
only to show the fact that
$(\mathcal{F}, \mathcal{E}_1)$
is compactly embedded into $L^{2p}(E;m)$
in the proof of Proposition 3.2 of \cite{MR4075013}.
Hence, by Corollary \ref{Cor_cptemb}
we may drop the assumption (A4) in \cite{MR4075013}.

\end{Rmk}

\vspace{1eM}

\begin{Ex}
[Killed Brownian motion in a domain $D$]
%\ \par
Suppose
$D\subset \mathbb{R}^d$ be
a domain with a smooth boundary satisfying
\begin{equation*}
\label{eq_killedBM}
 \lim_{x\in D, |x|\rightarrow \infty}
 m
 \bigl(D\cap B(x, 1)\bigr)
=
 0
,
\end{equation*}
where
$m$ is the Lebesgue measure on $\mathbb{R}^d$.
Let
$\partial$ be a point added to $D$ so that
$
 D_\partial
:=
 D\cup \{\partial\}
$
is the one-point compactification of $D$.
A killed Brownian motion $X$ in $D$ is the process
given by
\begin{equation*}
 X_t
=
 \begin{cases}
  B_t,      & t< \tau_D,
 \\
  \partial, & t\geq \tau_D,
 \end{cases}
\end{equation*}
where
$B$ is a Brownian motion on $\mathbb{R}^d$
and
$\tau_D = \inf\{t>0: B_t \not\in D\}$
is the exit time of $B$ from $D$.
Its Dirichlet form is
$(\frac{1}{2}\boldsymbol{D}, H^1_0(D))$,
where
$H^1_0(D)$ is the Sobolev space
with zero boundary values.
It is known that
$X$ satisfies the resolvent strong Feller property
and the tightness property
(see \cite[Lemma 3.3]{MR3685590} for example).

Hence,
by combining Example \ref{Ex_BM_pdelKato}
with Corollary \ref{Cor_cptemb},
we can see that
$H^1_0(D)$ is compactly embedded into $L^{2p}(D)$
for $p\in [1, \infty)$ with
$d-p(d-2)>0$.
This is exactly
the classical Rellich-Kondrachov embedding theorem.
\end{Ex}

%\newpage
%\vspace{1eM}

%%%%%%%%%%%%%%%%%%%%%%%%%%%%%%%%%%%%%%%%%%%%%%%%%%%%%%%
\subsection{Sobolev embedding implies $L^p$-Kato}
\label{Sec_Sob_pKato}
%%%%%%%%%%%%%%%%%%%%%%%%%%%%%%%%%%%%%%%%%%%%%%%%%%%%%%%

In this section,
we prove \eqref{eq_Sob_Dyn}
introduced in Section \ref{Sec_intro},
that is, for
a measure $\mu\in \mathcal{D}^1(X)$ and
$1\leq p< p'$,
if $(\mathcal{F}, \mathcal{E}_1)$ is continuously
embedded in $L^{2p'}(E;\mu)$,
then
$\mu\in \mathcal{D}^p(X)$.
Furthermore,
we give variants of \eqref{eq_Sob_Dyn}
which give a sufficient condition for
$\mu$ belonging to $\mathcal{K}^p(X)$ or
$\mathcal{K}^{p, \delta}(X)$ for a suitable $\delta$.

\vspace{1eM}

We first consider the case $\mu=m$ and
consider the assertion that
$\mathcal{E}_1$ is replaced by $\mathcal{E}$:

\begin{Lem}
\label{Lem_Sob_pKato_2}
%\ \par
Let $p'\in (1, \infty)$ and
suppose that the following Sobolev type inequality
holds: there exists a constant $S>0$ such that
\begin{equation}
\label{eq_2p'Sob_m_0}
 \|u\|_{L^{2p'}(E; m)} ^2
\leq
 S
 \mathcal{E}(u, u)
\quad
 \text{for all }u\in \mathcal{F}
.
\end{equation}

\noindent
Then $m\in \mathcal{K}^{p, \delta}(X)$
for any $p\in [1, p')$ with
$\delta = 1 -\frac{p'}{p'-1} \frac{p-1}{p}$.
\end{Lem}

\begin{proof}
%\ \par
First,
by \cite{MR803094}, \eqref{eq_2p'Sob_m_0}
implies the ultra-contractivity,
that is, there exists $C>0$ such that
\begin{equation*}
 \|P_t\|_{L^1\rightarrow L^\infty}
\leq
 C
 t^{-\frac{p'}{p'-1}}
\quad
 \text{for all }t>0
,
\end{equation*}
where $\|\cdot\|_{L^q\rightarrow L^r}$
is the operator norm
from $L^q(E; m)$ to $L^r(E; m)$.
By Jensen's inequality, we then have
\begin{equation*}
 \|P_t\|_{L^{\frac{p}{p-1}}\rightarrow L^\infty}
\leq
 C^{\frac{p-1}{p}}
 t^{-\frac{p'}{p'-1}\frac{p-1}{p}}
\quad
 \text{for all }t>0
,
\end{equation*}
when $p\in (1, p')$.
Since
the transition kernel $P_t$ is $L^\infty$-contractive,
the above inequality also holds when $p=1$,
where we use the convention $1/0 = \infty$.

Fix $\alpha > 0$.
For any non-negative Borel function $f$ with
$f\in L^1(E; m)\cap L^\infty(E; m)$,
we have

\noindent
\begin{align*}
 R_\alpha f(x)
=&
 \int_0 ^\infty
  e^{-\alpha t}
  P_t f(x)
 dt
\\
\leq&
 \int_0 ^\infty
  e^{-\alpha t}
  \|P_t f \|_{L^\infty(E; m)}
 dt
\\
\leq&
 C^{\frac{p-1}{p}}
 \|f\|_{L^{\frac{p}{p-1}}(E; m)}
 \int_0 ^\infty
  e^{-\alpha t}
  t^{-\frac{p'}{p'-1}\frac{p-1}{p}}
 dt
\\
=&
 C^{\frac{p-1}{p}}
 \|f\|_{L^{\frac{p}{p-1}}(E; m)}
 \Gamma(\delta)
 \alpha^{-\delta}
\end{align*}

\noindent
for $m$-a.e. $x\in E$,
where
$\Gamma$ is the Gamma function
and
$
 \delta
=
 1-\frac{p'}{p'-1}\frac{p-1}{p}
>
 0
$.
Since
$
 R_\alpha f
=
 \int_E
  r_\alpha (\cdot, y) f(y)
 m(dy)
$
is $\alpha$-excessive
and the absolute continuity condition \eqref{eq_AC}
holds,
we have for every $x\in E$,

\noindent
\begin{align}
\notag
 R_\alpha f(x)
=&
 \lim_{\varepsilon \downarrow 0}
 e^{-\alpha \varepsilon}
 \mathbb{E}_x [R_\alpha f (X_\varepsilon)]
=
 \lim_{\varepsilon \downarrow 0}
 e^{-\alpha \varepsilon}
 \int_E
  p_\varepsilon(x, y)
  R_\alpha f (y)
 m(dy)
\\
\label{eq_excessive}
\leq&
 C^{\frac{p-1}{p}}
 \|f\|_{L^{\frac{p}{p-1}}(E; m)}
 \Gamma(\delta)
 \alpha^{-\delta}
.
\end{align}

Now,
fix $x\in E$, $M>0$ and a compact set $K\subset E$.
By applying \eqref{eq_excessive} with
$
 f
=
 1_K(\cdot)
 \bigl(
  r_\alpha(x, \cdot)\wedge M
 \bigr)^{p-1}
$,
we have

\noindent
\begin{align*}
 \int_K
  \bigl(
   r_\alpha(x, y)\wedge M
  \bigr)^p
 m(dy)
\leq&
 \int_K
  \bigl(
   r_\alpha(x, y)\wedge M
  \bigr)^{p-1}
  r_\alpha(x, y)
 m(dy)
\\
\leq&
 C^{\frac{p-1}{p}}
 \biggl(
  \int_K
   \bigl(
    r_\alpha(x, y)\wedge M
   \bigr)^{p}
  m(dy)
 \biggr)^{\frac{p-1}{p}}
 \Gamma(\delta)
 \alpha^{-\delta}
,
\end{align*}

\noindent
which means that
\begin{equation*}
 \biggr(
  \int_K
   \bigl(
    r_\alpha(x, y)\wedge M
   \bigr)^{p}
  m(dy)
 \biggl)^{\frac{1}{p}}
\leq
 C^{\frac{p-1}{p}}
 \Gamma(\delta)
 \alpha^{-\delta}
.
\end{equation*}

\noindent
Let $M\uparrow \infty$ and $K\uparrow E$.
Hence
we have the conclusion
$m\in \mathcal{K}^{p, \delta}(X)$
by the dominated convergence theorem.
\end{proof}

\vspace{1eM}

%%%%%%%%%%%%%%%%%%%%%%%%%%%%%%%%%%%%%%%%%%%%%%%%%%%%%%%

The following is our second main theorem,
which states the opposite of Theorem \ref{Thm_Sobolev}
in some sense.

\begin{Thm}
\label{Thm_Sob_pKato}
%\ \par
Let $p'\in (1, \infty)$
and
let $\mu$ be a measure in $\mathcal{D}^1(X)$.
%whose fine support is identical
%to the topological support.

\vspace{-3mm}

\begin{itemize}
\setlength{\itemsep}{0mm}

\item[(i)]
Suppose the following
Sobolev type inequality holds:
there exists a constant $S>0$ such that

\noindent
\begin{equation}
\label{eq_2p'Sob_mu}
 \|u\|_{L^{2p'}(E; \mu)} ^2
\leq
 S
 \mathcal{E}_1(u, u)
\end{equation}
for all $u\in \mathcal{F}$.
Then
$\mu\in \mathcal{D}^{p}(X)$ for any $p\in [1, p')$.

\item[(ii)]
Suppose the following
Sobolev type inequality holds:
there exists a function
$K:(0, \infty) \rightarrow (0, \infty)$
with
$\eps^{-1} K(\eps) \uparrow \infty$
as $\eps \downarrow 0$,
such that

\noindent
\begin{equation}
\label{eq_2p'_Kato_Sob}
 \|
  u
 \|_{L^{2p'}(E;\mu)} ^2
\leq
 \eps
 \mathcal{E}_1(u, u)
+
 K(\eps)
 \|
  u
 \|_{L^{2}(E;m)} ^2
\end{equation}
for any $u\in \mathcal{F}$ and $\eps > 0$.
Then
$\mu\in \mathcal{K}^{p}(X)$ for any $p\in (1, p')$.

\item[(iii)]
Suppose the following
Sobolev type inequality holds:
there exist constants $A>0$ and $\theta\in (0, 1]$
such that

\noindent
\begin{equation}
\label{eq_2p'_thetaKato_Sob}
 \|
  u
 \|_{L^{2p'}(E;\mu)}
\leq
 A
 \sqrt{\mathcal{E}_1(u, u)}^{(1-\theta)}
 \|
  u
 \|_{L^{2}(E;m)} ^{\theta}
\end{equation}
for any $u\in \mathcal{F}$.
Then
$\mu\in \mathcal{K}^{p, \theta(1-\delta)}(X)$
for any $p\in (1, p')$ with
$\delta = 1 -\frac{p'}{p'-1} \frac{p-1}{p}$.

\end{itemize}
\end{Thm}

%%%%%%%%%%%%%%%%%%%%%%%%%%%%%%%%%%%%%%%%%%%%%%%%%%%%%%%

\begin{Rmk}
\label{Rmk_Sob_pKato}
\ \par
\vspace{-3mm}

\begin{enumerate}
\setlength{\itemsep}{0mm}

\item
In claim (i),
when $p=p'$ the result does not hold in general.
See Example \ref{Ex_BM_Sobolev} below.

\item
Regarding claims (ii) and (iii),
our proofs do not work when $p=1$ and $p=p'$,
but
there are examples that support
that the conclusion also holds for such $p$.
Indeed,
in \cite[Theorem 4.9]{MR644024},
Aizenman and Simon proved that
for a Brownian motion on $\mathbb{R}^d$,
the inequality of the type \eqref{eq_2p'_Kato_Sob}
with $p'=1$ and
$
 K(\eps)
=
 A
 \exp
 \bigl\{
  e^{B \eps^{-a}}
 \bigr\}
$
implies that
$\mu$ belongs to $\mathcal{K}^1(X)$.
For the case $1<p=p'$,
see Example \ref{Ex_BM_Sobolev} below.

\item
Regarding claim (iii),
it can be seen that
the exponent $\theta(1-\delta)$ is appropriate
by the following direct calculations:
Let $q>1$, $\rho\in (0, 1]$ and
$\mu\in \mathcal{K}^{q, \rho}(X)$.
For $1<p<p'<q$,
set
$\theta' = \frac{q}{q-1}\frac{p'-1}{p'}$
and
$\delta = 1-\frac{p'}{p'-1}\frac{p-1}{p}$.
Note that
$\frac{1}{p'} = 1-\theta' + \frac{\theta'}{q}$
and
$\frac{1}{p} = \delta + \frac{1-\delta}{p'}$.
We also set $\theta = \rho\theta'$.
Then,
H\"older's inequality and
Theorem \ref{Thm_Sobolev} (i) and (iii) imply that

\noindent
\begin{align*}
 \|
  u
 \|_{L^{2p'}(E;\mu)}
\leq&
 \|
  u
 \|_{L^{2}(E;\mu)} ^{1-\theta'}
 \|
  u
 \|_{L^{2q}(E;\mu)} ^{\theta'}
\\
\leq&
 \Bigl(
  \gamma(1)
  \sqrt{\mathcal{E}_1(u, u)}
 \Bigr)^{1-\theta'}
 \Bigl(
  B
  \sqrt{\mathcal{E}_1(u, u)}^{1-\rho}
  \|
   u
  \|_{L^{2}(E;m)} ^\rho
 \Bigr)^{\theta'}
\\
=&
 \gamma(1)^{1-\theta'}
 B^{\theta'}
 \sqrt{\mathcal{E}_1(u, u)}^{1-\theta}
 \|
  u
 \|_{L^{2}(E;m)} ^{\theta}
,
\end{align*}
that is, the assumption (4.9) holds.
On the other hand,
H\"older's inequality also implies that

\noindent
\begin{align*}
&
 \biggl(
  \sup_{x\in E}
  \int_E
   r_\alpha(x, y)^p
  \mu(dy)
 \biggr)^{\frac{1}{p}}
\\
\leq&
 \biggl(
  \sup_{x\in E}
  \int_E
   r_\alpha(x, y)
  \mu(dy)
 \biggr)^{\delta}
 \biggl(
  \sup_{x\in E}
  \int_E
   r_\alpha(x, y)^{p'}
  \mu(dy)
 \biggr)^{\frac{1-\delta}{p'}}
\\
\leq&
 \biggl(
  \sup_{x\in E}
  \int_E
   r_\alpha(x, y)
  \mu(dy)
 \biggr)^{\delta+(1-\theta')(1-\delta)}
 \biggl(
  \sup_{x\in E}
  \int_E
   r_\alpha(x, y)^{q}
  \mu(dy)
 \biggr)^{\frac{\theta'(1-\delta)}{q}}
\\
=&
 O(\alpha^{-\rho\theta'(1-\delta)})
\end{align*}
as $\alpha\rightarrow \infty$,
that is,
the conclusion
$\mu\in \mathcal{K}^{p, \theta(1-\delta)}(X)$
holds.

\end{enumerate}
\end{Rmk}

%\vspace{1eM}

%%%%%%%%%%%%%%%%%%%%%%%%%%%%%%%%%%%%%%%%%%%%%%%%%%%%%%%

\begin{proof}
[Proof of Theorem \ref{Thm_Sob_pKato}]
%\ \par
We first prove (i).
Let $X^{(1)}$ be the 1-subprocess of $X$.
As we have seen
in the proof of Proposition \ref{Prop_timechange},
it holds that $\mu\in S_1(X^{(1)})$.
%and its fine support with respect
%to $Y$ is identical to the topological support.
We also have
$\mu \in \mathcal{D}^1(X^{(1)})$.

Denote $Y$ the time changed process of $X^{(1)}$
with respect to the PCAF of $X^{(1)}$
with Revuz measure $\mu$
and denote
its Dirichlet form as
$(\mathcal{E}^Y, \mathcal{F}^Y)$.
Since
the Dirichlet form of $X^{(1)}$ is
$(\mathcal{E}_1, \mathcal{F})$,
\eqref{eq_2p'Sob_mu} implies that
\begin{equation*}
 \|u\|_{L^{2p'}(\mu)} ^2
\leq
 S
 \mathcal{E}^Y(u, u)
\quad
 \text{for all }u\in \mathcal{F}^Y
.
\end{equation*}

\noindent
Let $p\in [1, p')$.
By
Lemma \ref{Lem_Sob_pKato_2}
we have $\mu\in \mathcal{D}^{p}(Y)$,
and then by
Proposition \ref{Prop_timechange} we have
$\mu\in \mathcal{D}^{p}(X^{(1)})$.
The conclusion
follows from
Remark \ref{Rmk_subproc}.

\vspace{1eM}

We next prove (ii).
Let $\beta > 0$ and
denote
$Y$ the time changed process of
the $\beta$-subprocess $X^{(\beta)}$
with respect to the PCAF of $X^{(\beta)}$
with Revuz measure $\mu$.
The $\alpha$-order resolvent of $Y$ can be written as
\begin{equation*}
 R^\beta_{\alpha, A}f(x)
=
 \mathbb{E}_x
 \biggl[
  \int_0 ^\infty
   e^{-\alpha A^\mu_t -\beta t}
   f(X_t)
  dA^\mu_t
 \biggr]
,
\quad
 \alpha>0, x\in F
.
\end{equation*}

\noindent
By Lemma \ref{Lem_AC}, $Y$ satisfies the absolute
continuity condition \eqref{eq_AC}.
Denote
the $\alpha$-order resolvent kernel of $Y$ as
$r^\beta_\alpha(x, y)$.
Note that
for all $x\in F$ and $\alpha, \beta>0$,
it holds that
\begin{equation*}
 r_\beta(x, y)
=
 r^\beta_\alpha(x, y)
+
 \alpha
 \int_F
  r_\beta(x, z)
  r^\beta_\alpha(z, y)
 \mu(dz)
\quad
 \text{ for $\mu$-a.e. }
 y\in F
\end{equation*}
(see \cite[(6.4.11)]{MR2778606} for example).
By an argument similar to that leading to \eqref{eq_pResol},
we have
\begin{align}
\label{eq_beta_bdd}
\begin{split}
&
 \biggl(
  \sup_{x\in F}
  \int_E
   r_\beta(x, y)^p
  \mu(dy)
 \biggr)^{\frac{1}{p}}
\\
\leq&
 \biggl(
  1
 +
  \alpha
  \sup_{x\in E}
  \int_E
   r_\beta(x, y)
  \mu(dy)
 \biggr)
 \biggl(
  \sup_{z\in F}
  \int_F
   r^\beta_\alpha(z, y)^p
  \mu(dy)
 \biggr)^{\frac{1}{p}}
.
\end{split}
\end{align}

\noindent
By a similar argument as to obtain \eqref{eq_EbddF},
we also have
\begin{equation}
\label{eq_beta_EbddF}
 \sup_{x\in E}
 \int_E
  r_\beta(x, y)^p
 \mu(dy)
=
 \sup_{x\in F}
 \int_E
  r_\beta(x, y)^p
 \mu(dy)
.
\end{equation}

\noindent
Define the process $Z$ by $Z_t := Y_{\eps^{p'} t}$.
$Z$ is an $\eps^{-p'}\mu$-symmetric Hunt process
on $F$,
and its $\alpha$-order resolvent kernel
can be written as
$r^\beta_{\eps^{-p'}\alpha}(x, y)$.
Denote by
$(\mathcal{E}^Z, \mathcal{F}^Z)$
the associated Dirichlet form on $L^2(F; \eps^{-p'}\mu)$.

Now,
by setting $\beta = \eps^{-1}K(\eps)$,
\eqref{eq_2p'_Kato_Sob} implies that
\begin{equation*}
 \|
  u
 \|_{L^{2p'}(E;\eps^{-p'}\mu)} ^2
\leq
 \mathcal{E}^Z(u, u)
\quad
 \text{for all }
 u\in \mathcal{F}^Z
.
\end{equation*}

\noindent
By the same argument as the proof of
Lemma \ref{Lem_Sob_pKato_2},
there exists a constant $C>0$,
which is independent of $\eps$ and $\beta$
by virtue of \cite{MR803094},
such that
\begin{equation*}
 \biggl(
  \sup_{x\in F}
  \int_F
   r^\beta_{\eps^{-p'}\alpha} (x, y)^p
  \eps^{-p'}
  \mu(dy)
 \biggr)^{\frac{1}{p}}
\leq
 C^{\frac{p-1}{p}}
 \Gamma(\delta)
 \alpha^{-\delta}
\end{equation*}

\noindent
for all $\alpha > 0$.
In particular,
by substituting $\alpha = \eps^{p'}$, we have
\begin{equation*}
 \biggl(
  \sup_{x\in F}
  \int_F
   r^\beta_{1} (x, y)^p
  \mu(dy)
 \biggr)^{\frac{1}{p}}
\leq
 C^{\frac{p-1}{p}}
 \Gamma(\delta)
 \eps^{\frac{p'}{p}-p'\delta}
=
 C^{\frac{p-1}{p}}
 \Gamma(\delta)
 \eps^{1-\delta}
,
\end{equation*}
where we used the relation
$
 \frac{1}{p}
=
 \delta
+
 \frac{1-\delta}{p'}
$.

Take $\eps>0$
sufficiently small so that $\beta>1$.
By
combining the above inequality with
\eqref{eq_beta_bdd}, \eqref{eq_beta_EbddF}
and
the monotonicity $r_\beta(x, y)\leq r_1(x, y)$,
we have
\begin{align}
\label{eq_beta_bdd_C}
 \biggl(
  \sup_{x\in E}
  \int_E
   r_\beta(x, y)^p
  \mu(dy)
 \biggr)^{\frac{1}{p}}
\leq
 \biggl(
  1
 +
  \sup_{x\in E}
  \int_E
   r_1(x, y)
  \mu(dy)
 \biggr)
 C^{\frac{p-1}{p}}
 \Gamma(\delta)
 \eps^{1-\delta}
,
\end{align}
hence
we obtain the desired conclusion
$\mu\in \mathcal{K}^p(X)$.

\vspace{1eM}

We finally prove (iii).
By the assumption \eqref{eq_2p'_thetaKato_Sob},
It holds that the inequality \eqref{eq_2p'_Kato_Sob}
with
$
 K(\eps)
=
 A^{\frac{2}{\theta}}
 \eps^{-\frac{1-\theta}{\theta}}
$.
By setting
$
 \beta
=
 \eps^{-1}K(\eps)
=
 A^{\frac{2}{\theta}}
 \eps^{-\frac{1}{\theta}}
$,
we have
$
 \eps
=
 A^{2}
 \beta^{-\theta}
$
and then by \eqref{eq_beta_bdd_C}
\begin{align*}
&
 \biggl(
  \sup_{x\in E}
  \int_E
   r_\beta(x, y)^p
  \mu(dy)
 \biggr)^{\frac{1}{p}}
\\
\leq&
 \biggl(
  1
 +
  \sup_{x\in E}
  \int_E
   r_1(x, y)
  \mu(dy)
 \biggr)
 C^{\frac{p-1}{p}}
 \Gamma(\delta)
 A^{2(1-\delta)}
 \beta^{-\theta(1-\delta)}
.
\end{align*}

\noindent
Therefore we have the desired conclusion
$\mu\in \mathcal{K}^{p, \theta(1-\delta)}(X)$.
\end{proof}

\vspace{1eM}

%%%%%%%%%%%%%%%%%%%%%%%%%%%%%%%%%%%%%%%%%%%%%%%%%%%%%%%

\begin{Ex}
[Brownian motion]
\label{Ex_BM_Sobolev}
We continue with Example \ref{Ex_BM_pdelKato}.
As in Section \ref{Sec_intro},
the classical Sobolev embedding theorem
on $\mathbb{R}^d$ gives that
$H^1(\mathbb{R}^d)$ is continuously embedded
into $L^{2p}(\mathbb{R}^d)$
for
$p\in [1, \infty)$ with $d-p(d-2)\geq 0$.
By combining this with Theorem \ref{Thm_Sob_pKato} (i),
the Lebesgue measure on $\mathbb{R}^d$ is in
$\mathcal{D}^{p}(X)$
for
$p\in [1, \infty)$ with $d-p(d-2)>0$.

In particular,
when $d\geq 3$,
by setting $p^* = d/(d-2)$
the critical Sobolev embedding theorem gives that
$H^1(\mathbb{R}^d)$ is continuously embedded into
$L^{2p}(\mathbb{R}^d)$
for
$p\in [1, p^*)$,
but
the Lebesgue measure on $\mathbb{R}^d$
does not belong to $\mathcal{D}^{p^*}(X)$.
This is an example of Remark \ref{Rmk_Sob_pKato}.1.

We now assume $d\geq 1$.
We also assume that $d-p'(d-2)\geq 0$ and $1<p<p'$.
The classical
Gagliardo-Nirenberg interpolation inequality
says that, by setting
$
 \theta
=
 1-\frac{d(p'-1)}{2p'}
=
 \frac{d - p'(d-2)}{2p'}
$,
there exists a positive constant $C$ such that

\noindent
\begin{align*}
 \|
  u
 \|_{L^{2p'}(\mathbb{R}^d)}
\leq
 C
 \|
  \nabla u
%\|_{L^{2}(\mathbb{R}^d)} ^{\frac{d(p'-1)}{2p'}}
 \|_{L^{2}(\mathbb{R}^d)} ^{1-\theta}
 \|
  u
%\|_{L^{2}(\mathbb{R}^d)} ^{1-\frac{d(p'-1)}{2p'}}
 \|_{L^{2}(\mathbb{R}^d)} ^{\theta}
\end{align*}
holds for all $u\in H^{1}(\mathbb{R}^d)$.
By applying
Theorem \ref{Thm_Sob_pKato} (iii)
and then by letting $p'\downarrow p$,
we can see that
the Lebesgue measure on $\mathbb{R}^d$
belongs to
$
 \mathcal{K}^{p, \frac{d-p(d-2)}{2p}-\varepsilon}(X)
$
for all
$\varepsilon \in (0, \frac{d-p(d-2)}{2p})$,
which is
a slightly weaker condition than
the conclusion of Example \ref{Ex_BM_pdelKato}.
By
combining with the Gagliardo-Nirenberg inequality,
we can also see that
the fact that
the Lebesgue measure on $\mathbb{R}^d$ belongs to
$
 \mathcal{K}^{p, \frac{d-p(d-2)}{2p}}(X)
$,
which
is an example of Remark \ref{Rmk_Sob_pKato}.2
for $1<p=p'$.
\end{Ex}

%\newpage
\vspace{1eM}

%%%%%%%%%%%%%%%%%%%%%%%%%%%%%%%%%%%%%%%%%%%%%%%%%%%%%%%
\section{Application: Continuity of the intersection measure  in time}
\label{Sec_ISmeas}
%%%%%%%%%%%%%%%%%%%%%%%%%%%%%%%%%%%%%%%%%%%%%%%%%%%%%%%

In this section,
we give an application of the $L^p$-Kato class
to the continuity of the intersection measure in time.
Throughout this section,
we assume that $p\geq 2$ is an integer
and
that the reference measure $m$ is
in the $L^p$-Dynkin class.
Let
$X^{(1)}, \ldots, X^{(p)}$ be
independent Hunt processes with the same law as $X$.
We write
$\zeta^{(1)},\ldots, \zeta^{(p)}$
as their life times
and write
$x^{(1)}_0,\ldots, x^{(p)}_0$
as their starting points, respectively.

First,
we review the construction of the intersection measure.
For detail, see
\cite{MR2584458, MR4075013} for example.
Fix
bounded Borel sets
$
 J^{(1)},\ldots, J^{(p)}
\subset
 [0, \infty)
$
and write
$
 J
=
 \prod_{i=1} ^p
 J^{(i)}
$.
For each $\varepsilon >0$,
we define
{\it the approximated (mutual) intersection measure}
$
 \ell^{\mathrm{IS}}_{J, \varepsilon}
$
of $X^{(1)}, \ldots, X^{(p)}$
with respect to
the (multi-parameter) time interval $J$
by
\begin{equation*}
 \langle
  \ell^{\mathrm{IS}}_{J, \varepsilon}
 ,
  f
 \rangle
=
 \int_E
  f(x)
  \biggl[
   \prod_{i=1}^p
   \int_{J^{(i)}}
    p_\varepsilon(x, X^{(i)}_s)
   ds
  \biggr]
 m(dy)
\end{equation*}
for $f\in \mathcal{B}_b(E)$,
where,
for convenience we regard
$
 p_\varepsilon(x, X^{(i)}_s)
=
 0
$
when $s\geq \zeta^{(i)}$.
Then,
there exists a random measure
$\ell^{\mathrm{IS}}_J$ on $E$
such that,
$
 \ell^{\mathrm{IS}}_{J, \varepsilon}
$
converges vaguely to
$\ell^{\mathrm{IS}}_J$ in $\mathcal{M}(E)$
and that
\begin{equation*}
 \lim_{\varepsilon \rightarrow 0}
 \mathbb{E}
 \bigl[
  |
   \langle
    f, \ell^{\mathrm{IS}}_{J, \varepsilon}
   \rangle
  -
   \langle
    f, \ell^{\mathrm{IS}}_{J}
   \rangle
  |^k
 \bigr]
=
 0
\end{equation*}

\noindent
for any integer $k\geq 1$ and $f\in C_0(E)$,
where
$\mathcal{M}(E)$ is the set of Radon measures on $E$
equipped with the vague topology.
We call
the limit $\ell^{\mathrm{IS}}_J$ as
{\it
 the (mutual) intersection measure
 of $X^{(1)}, \ldots, X^{(p)}$
 with respect to $J$%
}.
For
$
 \boldsymbol{t}
=
 (t_1, \ldots, t_p)
\in
 [0, \infty)^p
$
we simply denote the approximated intersection measure
and the intersection measure
with respect to
$
 [\boldsymbol{0},\boldsymbol{t}]
:=
 \prod_{i=1} ^p
 [0, t_i]
$
as
$
 \ell^{\mathrm{IS}}_{\boldsymbol{t}, \varepsilon}
$
and
$
 \ell^{\mathrm{IS}}_{\boldsymbol{t}}
$,
respectively.

The intersection measure
$
 \ell^{\mathrm{IS}}_{J}
$
enjoys the so-called Le Gall's moment formula:
for any $f\in \mathcal{B}_b(E)$ with compact support
and
for any integer $k\geq 1$,
it holds that

\noindent
\begin{align}
\notag
&
 \mathbb{E}
 \bigl[
  \langle
   f
  ,
   \ell^{\mathrm{IS}} _{J}
  \rangle^k
 \bigr]
\\
\begin{split}
=&
 \int_{E^k}
  f(x_1)\cdots f(x_k)
  \prod_{i=1}^p
  \biggl\{
   \sum_{\sigma\in \mathfrak{S}_k}
   \int_{(J^{(i)})^k _<}
    \prod_{j=1} ^k
    p_{s_j - s_{j-1}}
    (
     x_{\sigma(j-1)}
    ,
     x_{\sigma(j)}
    )
\\
&\hspace{60mm}
   ds_1 \cdots ds_k
  \biggr\}
 m(dx_1)\cdots m(dx_k)
,
\end{split}
\end{align}
where
$
 (J^{(i)})^k _<
:=
 \{
  (s_1, \ldots, s_k)\in (J^{(i)})^k
 ;
  s_1 < \cdots < s_k
 \}
$
and
$\mathfrak{S}_k$ is the set of permutations of
$\{1, \ldots, k\}$.
For convenience we set
$\sigma(0) = 0$ for $\sigma\in \mathfrak{S}_k$
and set
$x_0 = x^{(i)}_0$.

\vspace{1eM}

The goal of this section is to realize
$
 \{
  \ell^{\mathrm{IS}}_{\boldsymbol{t}}
 :
  \boldsymbol{t}\in [0,\infty)^p
 \}
$
as a measure-valued continuous stochastic process.
That is,

\begin{Thm}
\label{Thm_ISconti}
%\ \par
Assume
$m \in \mathcal{K}^{p, \delta}(X)$
for some $\delta\in (0, 1]$.
Then it holds that
\vspace{-3mm}

\begin{itemize}
\setlength{\itemsep}{0mm}

\item[(i)]
the $\mathcal{M}(E)$-valued process
$
 \{
  \ell^{\mathrm{IS}}_{\boldsymbol{t}}
 :
  \boldsymbol{t}\in [0,\infty)^p
 \}
$
has a continuous modification,

\item[(ii)]
for any $f\in \mathcal{B}_b(E)$,
the real-valued process
$
 \{
  \langle
   f
  ,
   \ell^{\mathrm{IS}} _{\boldsymbol{t}}
  \rangle
 :
  \boldsymbol{t}\in [0, \infty)^p
 \}
$
has a modification
whose paths are locally
$\gamma$-H\"older continuous of every order
$\gamma \in (0, \delta)$.
\end{itemize}
\end{Thm}

\begin{proof}
%\ \par
Our proof is based on
that of \cite[Lemma 2.2.4]{MR2584458}.
First, we claim the following estimate:
for any integer $k\geq 1$,
positive constant $T>0$ and
$
 \boldsymbol{s}
,
 \boldsymbol{t}
\in
 [0, T]^p
$,
it holds that

\noindent
\begin{align}
\label{eq_boldst^k}
 \mathbb{E}
 \bigl[
  |
  \langle
   f
  ,
   \ell^{\mathrm{IS}}_{\boldsymbol{t}}
  \rangle
 -
  \langle
   f
  ,
   \ell^{\mathrm{IS}}_{\boldsymbol{s}}
  \rangle
  |^k
 \bigr]
\leq&
 (k!)^p
 \Bigl(
  2^{p}
  \|f\|_\infty
  (\eta(T)+1)^{p}
  \sup_{0<t\leq pT}
  \{
   t^{-\delta}\eta(t)
  \}
 \Bigr)^k
 |\boldsymbol{t}-\boldsymbol{s}|^{\delta k}
,
\end{align}
where
$\eta(T)$ is introduced in \eqref{eq_eta}
and
$
 |\boldsymbol{t}-\boldsymbol{s}|
$
is the Euclidean distance between
$\boldsymbol{t}$ and $\boldsymbol{s}$
in $[0, \infty)^p$.

Fix
$
 \boldsymbol{s}
,
 \boldsymbol{t}
\in
 [0, T]^p
$.
Define the family
$
 \{
  J_l
 \}_{l=1} ^{2^p-1}
$
consisting of products of intervals as

\noindent
\begin{align*}
&
 \{
  J_l
 \}_{l=1} ^{2^p-1}
\\
=&
 \biggl\{
  \prod_{i=1} ^p
  [a_i, b_i]
 \hspace{1mm}
 \biggl|
 \hspace{1mm}
  [a_i, b_i]
 =
  [0            , s_i \wedge t_i]
 \text{ or }
  [s_i\wedge t_i, s_i \vee   t_i]
 ,
  1\leq i\leq p
 \biggr\}
\setminus
 \biggl\{
  \prod_{i=1} ^p
  [0          , s_i \wedge t_i]
 \biggr\}
.
\end{align*}
Then we have
$
 [\boldsymbol{0}, \boldsymbol{t}]
\triangle
 [\boldsymbol{0}, \boldsymbol{s}]
\subset
 \bigcup_l
 J_l
$
and then
\begin{align*}
 |
 \langle
  f
 ,
  \ell^{\mathrm{IS}}_{\boldsymbol{t}}
 \rangle
-
 \langle
  f
 ,
  \ell^{\mathrm{IS}}_{\boldsymbol{s}}
 \rangle
 |
\leq
 \sum_{l=1} ^{2^p-1}
 \langle
  |f|
 ,
  \ell^{\mathrm{IS}}_{J_l}
 \rangle
,
\quad
 \mathbb{P}\text{-a.s.}
\end{align*}

\noindent
By Le Gall's moment formula, we have
\begin{align*}
&
 \mathbb{E}
 \bigl[
  \langle
   |f|
  ,
   \ell^{\mathrm{IS}} _{\varepsilon, J_l}
  \rangle^k
 \bigr]
\\
=&
 \int_{E^k}
  |f|(x_1)\cdots |f|(x_k)
  \prod_{i=1}^p
  \biggl\{
   \sum_{\sigma\in \mathfrak{S}_k}
   H^{(i)}_{a_i, b_i}
   (x_{\sigma(1)}, \ldots, x_{\sigma(k)})
  \biggr\}
 m(dx_1)\cdots m(dx_k)
,
\end{align*}

\noindent
where
\begin{align*}
&
 H^{(i)}_{a_i, b_i}
 (x_1, \ldots, x_p)
\\
:=&
 \int_E
 \int_{[0, \infty)^p}
  1_{
   \bigl\{
    \sum_{j=1} ^k
    r_j
   \leq
    b_i - a_i
   \bigr\}
  }
  \prod_{j=1} ^k
  p_{r_j} (x_{j-1}, x_j)
 \hspace{1mm}
 dr_1 \cdots dr_k
 \hspace{1mm}
 \nu^{(i)}_{a_i}(dx_0)
\end{align*}
and
$
 \nu^{(i)}_{a_i}(dx_0)
=
 \mathbb{P}_{x^{(i)}_0}({X^{(i)}_{a_i} }\in dx_0)
$.
Note that
\begin{equation*}
 H^{(i)}_{a_i, b_i}
 (x_1, \ldots, x_p)
\leq
 \int_E
  \Bigl[
   \prod_{j=1} ^k
   \int_0^{b_i-a_i}
    p_{s} (x_{j-1}, x_j)
   ds
  \Bigr]
 \nu^{(i)}_{a_i}(dx_0)
\end{equation*}
and then
$
 \bigl\|
  H^{(i)}_{a_i, b_i}
 \bigr\|_{L^p(E^k, m^{\otimes k})}
\leq
 \eta(b_i - a_i)^k
$.
We have
by H\"older's inequality

\noindent
\begin{align*}
 \mathbb{E}
 \bigl[
  \langle
   |f|
  ,
   \ell^{\mathrm{IS}} _{\varepsilon, J_l}
  \rangle^k
 \bigr]
\leq&
 (k!)^p
 \|f\|_\infty ^k
 \prod_{i=1} ^p
 \bigl\|
  H^{(i)}_{a_i, b_i}
 \bigr\|_{L^p(E^k, m^{\otimes k})}
\\
\leq&
 (k!)^p
 \|f\|_\infty ^k
 \prod_{i=1} ^p
 \eta(b_i - a_i)^k
\leq
 (k!)^p
 \|f\|_\infty ^k
 \eta(T)^{(p-1)k}
 \eta(|\boldsymbol{t}-\boldsymbol{s}|)^k
.
\end{align*}
In the last inequality,
we used the fact that
$[a_i, b_i] = [s_i\wedge t_i, s_i \vee t_i]$
holds for at least one $i$
because of the definition of $J_l$.
Therefore

\noindent
\begin{align*}
 \mathbb{E}
 \bigl[
  |
  \langle
   f
  ,
   \ell^{\mathrm{IS}}_{\boldsymbol{t}}
  \rangle
 -
  \langle
   f
  ,
   \ell^{\mathrm{IS}}_{\boldsymbol{s}}
  \rangle
  |^k
 \bigr]
\leq&
 \biggl\{
 \sum_{l=1} ^{2^p-1}
 \mathbb{E}
 \bigl[
  \langle
   |f|
  ,
   \ell^{\mathrm{IS}}_{J_l}
  \rangle^k
 \bigr]^{1/k}
 \biggr\}^k
\\
\leq&
 (k!)^p
 \Bigl(
  2^p
  \|f\|_\infty
  (\eta(T)+1)^{p}
  \eta(|\boldsymbol{t}-\boldsymbol{s}|)
 \Bigr)^k
.
\end{align*}

\noindent
Taking account of Remark \ref{Rmk_pdelKato}.1,
the desired estimate \eqref{eq_boldst^k} follows from
$|\boldsymbol{t}-\boldsymbol{s}|\leq pT$
and the assumption $m\in\mathcal{K}^{p,\delta}(X)$.

\vspace{1eM}

Now,
$(\delta k - p)/k$ increases to $\delta$
as $k\rightarrow \infty$.
By applying
Kolmogorov's continuity theorem
(see \cite[Theorem 3.23]{MR1876169} for example)
to \eqref{eq_boldst^k},
we can find that the process
$
 [0, T]^p
\ni
 \boldsymbol{t}
\mapsto
 \langle
  f
 ,
  \ell^{\mathrm{IS}}_{\boldsymbol{t}}
 \rangle
\in
 \mathbb{R}
$
has a continuous modification
whose paths are $\gamma$-H\"older continuous of
every order $\gamma \in (0, \delta)$.
Since
$T>0$ is arbitrary, we obtain (ii).

To prove (i),
take a dense subset $\{\phi_n\}_{n=1} ^\infty$
in $C_0^+(E)$, the family of non-negative
continuous functions with compact support
equipped with the uniform metric.
Then
$\mathcal{M}(E)$ is homeomorphic to a subset of
$\mathbb{R}^\infty$ by the mapping
$
 \mathcal{M}(E)
\ni
 \mu
\mapsto
 \{
  \langle
   \mu, \phi_n
  \rangle
 \}_{n=1}^\infty
\in
 \mathbb{R}^\infty
$,
and hence the process
$
 [0, \infty)^p
\ni
 \boldsymbol{t}
\mapsto
 \ell^\mathrm{IS}_{\boldsymbol{t}}
\in
 \mathcal{M}(E)
$
has a continuous version by (ii).
Therefore we complete the proof.
\end{proof}

%%%%%%%%%%%%%%%%%%%%%%%%%%%%%%%%%%%%%%%%%%%%%%%%%%%%%%%
%%%Acknowledgements
%%%%%%%%%%%%%%%%%%%%%%%%%%%%%%%%%%%%%%%%%%%%%%%%%%%%%%%
\section*{Acknowledgements}
\addcontentsline{toc}{section}{Acknowledgements}
We
would like to thank
Professor Takashi Kumagai,
Professor Masayoshi Takeda
and an anonymous referee
for helpful discussions and comments.
We
also wish to thank Professor Kazuhiro Kuwae
for pointing out some mistakes in the first draft,
which finally improves our results.
This work
was supported by JSPS KAKENHI Grant Number JP18J21141.

%%%%%%%%%%%%%%%%%%%%%%%%%%%%%%%%%%%%%%%%%%%%%%%%%%%%%%%
%%%References
%%%%%%%%%%%%%%%%%%%%%%%%%%%%%%%%%%%%%%%%%%%%%%%%%%%%%%%
%\part*{}
\addcontentsline{toc}{section}{References}
%\bibliographystyle{abbrvalpha}
%\bibliographystyle{alpha,plain,jplain}

%%%%%%%%%%%%%%%%%%%%%%%%%%%%%%%%%%%%%%%%%%%%%%%%%%%%%%%
\end{document}